\documentclass[reqno,fleqn]{amsproc}%
\usepackage[frenchb,english]{babel}
\usepackage{amscd}
\usepackage{amsmath}
\usepackage{graphicx}
\usepackage{amsfonts}
\usepackage{amssymb}%
\theoremstyle{plain}
\newtheorem{theorem}{Theorem}[section]
\newtheorem{proposition}[theorem]{Proposition}
\newtheorem{lemma}[theorem]{Lemma}
\newtheorem{corollary}[theorem]{Corollary}
\newtheorem{definition}{Definition}

\theoremstyle{remark}

\newtheorem{exercise}{Exercise}[section]

\allowdisplaybreaks
\numberwithin{equation}{section}

\begin{document}
\title{Exceptional symmetric domains}
\author{Guy Roos}
\thanks{Lectures at the Workshop ``Several Complex Variables, Analysis on Complex Lie groups and Homogeneous
Spaces'', held at Zhejiang University, Hangzhou, China, Oct.~17-29, 2005. Revised version.}
\thanks{E-mail address: \texttt{guy.roos@normalesup.org }}

\begin{abstract}
We give the presentation of exceptional bounded symmetric domains using the
Albert algebra and exceptional Jordan triple systems.

The first chapter is devoted to Cayley-Graves algebras, the second to
exceptional Jordan triple systems. In the third chapter, we give a geometric
description of the two exceptional bounded symmetric domains, their
boundaries and their compactifications.
\end{abstract}

\maketitle
\tableofcontents


\section*{Introduction}



The classification of irreducible bounded symmetric complex domains is
well-known. They fall into four infinite series --- the ``classical
domains'' --- which can be defined as matrix spaces, using ordinary matrix
operations and classical linear groups, and two ``exceptional'' domains, of
respective complex dimension $16$ and $27$, which have no matrix description
(i.e., no description in a matrix space involving the usual matrix
operations).

The main purpose of these notes is to present an explicit algebraic and
geometric description of the two exceptional domains, which can no longer be
considered as ``unknown'', as well as some tools on them.

Analysis and geometry of classical domains have been extensively studied,
following the pioneer work of \index{Hua@\textsc{Hua, Luokeng} (1910--1985)}Hua Luokeng
\cite{Hua1963}, which consists of
a case-by-case study of the four classical series. A general theory for all
bounded complex domains also exists, using either semi-simple Lie groups
(see \cite{Helgason1978}, \cite{Koranyi1999}) or Jordan triple systems (see
\cite{Loos1977}, \cite{Roos1999}). The study of one particular classical
series still provides a good insight for conjecturing properties valid for
all bounded symmetric domains.

The explicit description of the exceptional domains, which was not known at
the time of Hua's book, has been available for at least 30 years. The
description involves $3\times3$ matrices with entries in the
\index{CayleyA@\textsc{Cayley, Arthur} (1821--1895)}%
\index{GravesJ@\textsc{Graves, John Thomas} (1806--1870)}Cayley--Graves
algebra $\mathbb{O}_{\mathbb{C}}$ of complex octonions. As this algebra is
\emph{non-associative}, these matrices do not carry the usual interpretation
of linear algebra theory and they do not build an associative matrix algebra
for the usual matrix operations. However, the space
\index{$H30@$\mathcal{H}_{3}(\mathbb{O}_{\mathbb{C}})$}%
$\mathcal{H}_{3}(\mathbb{%
O}_{\mathbb{C}})$ of such matrices which are Hermitian with respect to
Cayley conjugation can be endowed with the structure of a
\index{JordanP@\textsc{Jordan, Pascual} (1902--1980)}%
\index{Jordan algebra}\emph{Jordan
algebra}, using a product which generalizes in a natural way the symmetrized
product%
\begin{equation}
x\circ y=\frac{1}{2}(xy+yx)   \label{JordanProduct}
\end{equation}
of ordinary square matrices. This algebra is known as the
\index{AlbertAA@\textsc{Albert, Abraham Adrian} (1905--1972)}%
\index{Albert algebra}\emph{Albert
algebra} or \index{exceptional Jordan algebra}%
\emph{exceptional Jordan algebra}. It is the natural place to
describe the exceptional symmetric domain of dimension $27$. The second
exceptional symmetric domain (of complex dimension $16$) lives in the space $%
\mathcal{M}_{2,1}(\mathbb{O}_{\mathbb{C}})$ of $2\times1$ matrices with
octonion entries. This space has some analogy with the space $\mathcal{M}%
_{p,q}(\mathbb{C)}$ of ordinary rectangular matrices, endowed with the
\index{Jordan triple product}\emph{Jordan triple product}%
\begin{equation}
\left\{ xyz\right\} =xy^{\ast}z+zy^{\ast}x,   \label{JordanTripleProduct}
\end{equation}
where $y^{\ast}$ denotes the Hermitian adjoint (transposed conjugate) of $y$%
. The space \index{$M21O@$\mathcal{M}_{2,1}(\mathbb{O}_{\mathbb{C}})$}%
$\mathcal{M}_{2,1}(\mathbb{O}_{\mathbb{C}})$ also carries the
structure of a \index{Jordan triple system}\emph{Jordan triple system}, which allows an
algeb\-raico-geometric description of the exceptional domain of dimension $16
$.

The Jordan algebra $\mathcal{H}_{3}(\mathbb{O}_{\mathbb{C}})$ and the Jordan
triple system $\mathcal{M}_{2,1}(\mathbb{O}_{\mathbb{C}})$ are \emph{%
exceptional} not only because they are not part of an infinite series, but
more fundamentally because their algebraic products \emph{cannot} be related
to some associative product by formulas like (\ref{JordanProduct}) or (\ref%
{JordanTripleProduct}). But the explicit description of their algebraic
structure, combined with the general theory of Jordan triple systems and
bounded symmetric domains, provides easy access to the geometry and analysis
on the two exceptional symmetric domains. After this preliminary work, it
appears that exceptional domains are as easy (or not worse) to handle than
classical ones. It also appears that these two domains are as representative
as classical ones for exhibiting phenomena which lead to conjectures for all
symmetric domains.

\bigskip

The first chapter of these notes is devoted to Cayley--Graves algebras, the
second to exceptional Jordan triple systems. In the third chapter, we give a
geometric description of the two exceptional bounded symmetric domains,
their boundaries and their compactifications.


\section{Cayley algebras}



We denote by $k$ the field $\mathbb{R}$ of real numbers or the field $%
\mathbb{C}$ of complex numbers. A $k$\emph{-algebra} is a $k$-vector space $A
$, endowed with a $k$-bilinear \emph{product}
\begin{align*}
(x,y) & \mapsto xy \\
A\times A & \rightarrow A.
\end{align*}
This product is not assumed to be commutative nor associative. But we shall
assume it has a unit element $e\neq0$; this unit element is also denoted by $%
1$.

\subsection{Composition algebras}

\begin{definition}
\label{def1/1}A
\index{composition algebra}\emph{composition algebra}
(or \index{Hurwitz algebra}\index{HurwitzA@\textsc{Hurwitz, Adolf} (1859--1919)}%
\emph{Hurwitz algebra}) over $k$ is a pair $(A,n)$, where $A$ is a $k$-algebra and $n$ a
non-singular quadratic form on $A$, which is \emph{multiplicative} in the
sense that
\begin{equation}
n(ab)=n(a)n(b)\quad(a,b\in A).   \label{eq1/1}
\end{equation}
The form $n$ is called the \index{norm! on a composition algebra}\emph{norm}
of the composition algebra, and $n(a)$
is called the norm of $a$.
\end{definition}

It is clear that $n(e)=1$. We will identify $k$ and $ke$ using $\lambda
\mapsto\lambda e$. The elements of $ke$ are called the \emph{scalars} of $A$%
. For each $\lambda\in k$, we have $n(\lambda e)=\lambda^{2}$.

If $A=k$, there is a unique composition algebra structure over the $k$%
-vector space $k$, given by $n(\lambda)=\lambda^{2}$. This justifies the
above identification $\lambda\mapsto\lambda e$ in a general composition
algebra.

Denote by $A^{o}$ the opposite algebra of $A$ (i.e., the same vector space
with the opposite product $x\cdot y=yx$); clearly $(A^{o},n)$ is also a
composition algebra, which is called the \index{composition algebra!opposite ---}%
\emph{opposite composition algebra} to $(A,n)$.

In a composition algebra $(A,n)$, we denote by $(~:~)$ the bilinear form
associated to $n$ :
\begin{equation}
(a:b)=n(a+b)-n(a)-n(b)\quad(a,b\in A).   \label{eq1/2}
\end{equation}
Note there is no $\frac12$ factor in this definition, which implies $%
(a:a)=2n(a)$. Then the relation (\ref{eq1/1}) can also be written
\begin{equation*}
2(ab:ab)=(a:a)(b:b).
\end{equation*}
Polarizing this identity with respect to the variable $b$ yields
\begin{equation}
(ac:ad)=n(a)(c:d);   \label{eq1/3}
\end{equation}
polarizing with respect to the variable $a$ yields in the same way
\begin{equation}
(ac:bc)=(a:b)n(c).   \label{eq1/4}
\end{equation}
Polarizing again this last identity with respect to $c$, we obtain
\begin{equation}
(ac:bd)+(ad:bc)=(a:b)(c:d).   \label{eq1/5}
\end{equation}

Specializing this identity to $b\leftarrow1$, $c\leftarrow a$, and using (%
\ref{eq1/3}), we obtain for all $d\in A$
\begin{equation*}
(a^{2}:d)+n(a)(d:1)=(a:1)(a:d),
\end{equation*}
which is equivalent to
\begin{equation*}
(a^{2}-(a:1)a+n(a)1:d)=0.
\end{equation*}
As $n$ is assumed to be non-singular, this implies
\begin{equation*}
a^{2}-(a:1)a+n(a)1=0.
\end{equation*}
Define the \emph{trace} $t(a)$ in a composition algebra by
\begin{equation}
t(a)=(a:1).   \label{eq1/6}
\end{equation}
We have proved that each element $a$ in a composition algebra satisfies the
equation of degree $2$
\begin{equation}
a^{2}-t(a)a+n(a)1=0.   \label{eq1/7}
\end{equation}

\subsection{Cayley conjugation}

Let $(A,n)$ be a composition algebra. The \emph{(Cayley) conjugate} of an
element $a\in A$ is defined by
\begin{equation}
\widetilde{a}=(a:e)e-a.   \label{eq1/8}
\end{equation}
The
\index{Cayley conjugation}\emph{Cayley conjugation }$a\mapsto%
\widetilde{a}$ is the orthogonal symmetry (with respect to the quadratic
form $n$) which has $ke$ as its fixed point set. Therefore it is involutive
and isometric:
\begin{equation}
\left( \widetilde{a}\right)\widetilde{\ }=a,\quad n(\widetilde {a})=n(a).
\label{eq1/9}
\end{equation}
The defining relation (\ref{eq1/8}) can also be written $a+\widetilde{a}=t(a)
$; the identity (\ref{eq1/7}) is then equivalent to $-a\widetilde {a}+n(a)=0$%
. So norm, trace and conjugation are related by the relations
\begin{equation}
a+\widetilde{a}=t(a),\quad a\widetilde{a}=\widetilde{a}a=n(a)=n(\widetilde {a%
}).   \label{eq1/10}
\end{equation}
We also have, polarizing $n(\widetilde{a})=n(a)$,
\begin{equation}
(a:b)=(\widetilde{a}:\widetilde{b}).   \label{eq1/11}
\end{equation}

The relation (\ref{eq1/5}) with $b\leftarrow1$ gives
\begin{equation*}
(ac:d)+(ad:c)=((a:1)c:d);
\end{equation*}
as $(a:1)=a+\widetilde{a}$, we obtain $(ad:c)=(\widetilde{a}c:d)$. The
symmetric relation $(da:c)=(c\widetilde{a}:d)$ is proved in the same way.
These two identities can be better written as follows:
\begin{align}
(ax:y) & =(x:\widetilde{a}y),  \label{eq1/12} \\
(xa:y) & =(x:y\widetilde{a}).   \label{eq1/13}
\end{align}
Using these identities, we have $(ab:1)=(a:\widetilde{b})=(ba:1)$, that is
\begin{equation}
t(ab)=t(ba);   \label{eq1/14}
\end{equation}
we will say that the trace is ``commutative'' (with respect to the product).
Using again the identities (\ref{eq1/12})-(\ref{eq1/13}), we have
\begin{equation*}
t((ab)c)=(ab:\widetilde{c})=(a:\widetilde{c}\widetilde{b})=(ca:\widetilde {b}%
)=t((ca)b),
\end{equation*}
which means that $t((ab)c)$ is invariant under \emph{even} permutations of $%
(a,b,c)$. Using this fact and (\ref{eq1/14}), we get
\begin{equation*}
t((ab)c)=t((ca)b)=t((bc)a)=t(a(bc)),
\end{equation*}
that is, the trace is ``associative'' in the sense that
\begin{equation}
t((ab)c)=t(a(bc)).   \label{eq1/15}
\end{equation}
From (\ref{eq1/12})--(\ref{eq1/13}), we also have
\begin{equation*}
(\widetilde{ab}:c)=(ab:\widetilde{c})=(ca:\widetilde{b})=(c:\widetilde {b}%
\widetilde{a})
\end{equation*}
for all $c\in A$, which implies
\begin{equation}
(ab)\widetilde{\ }=\widetilde{b}\widetilde{a}.   \label{eq1/16}
\end{equation}
This means that the Cayley conjugation $a\mapsto\widetilde{a}$ is an
isomorphism from the composition algebra $(A,n)$ onto the opposite algebra $%
(A^{o},n)$.

Using (\ref{eq1/5}) and (\ref{eq1/12}), we get
\begin{equation*}
(a:b)(c:d)=(\widetilde{b}(ac):d)+(\widetilde{a}(bc):d)
\end{equation*}
for all $d\in A$, which implies, as $(~:~)$ is non-singular,
\begin{equation}
(a:b)c=\widetilde{b}(ac)+\widetilde{a}(bc).   \label{eq1/17}
\end{equation}
In the same way (or using the isomorphism with the opposite algebra), we
have
\begin{equation}
(a:b)c=(ca)\widetilde{b}+(cb)\widetilde{a.}   \label{eq1/18}
\end{equation}
Specializing these two relations to the case $a=b$, we get
\begin{equation}
n(a)c=\widetilde{a}(ac)=(ca)\widetilde{a}.   \label{eq1/19}
\end{equation}
The first equality can also be written $(\widetilde{a}a)c=\widetilde{a}(ac)$%
; using the fact that $a+\widetilde{a}$ is scalar, this implies $a^{2}c=a(ac)
$. One proves in the same way the identity $(ca)a=ca^{2}$. So we have proved
that the following identities:
\begin{equation}
a^{2}c=a(ac),\quad(ca)a=ca^{2}   \label{eq1/20}
\end{equation}
are verified in a composition algebra.

\begin{definition}
\label{def1/2}An algebra which satisfies the identities (\ref{eq1/20}) is
called an
\index{alternative algebra}\emph{alternative} algebra.
\end{definition}

The property of being alternative will be referred to as
\index{alternativity}\emph{alternativity}.

\subsection{Alternative algebras}

Let $A$ be a $k$-algebra. The \emph{commutator} $[x,y]$ and the \index{associator}\emph{%
associator} $[x,y,z]$ are respectively defined by
\begin{align*}
\lbrack x,y] & =xy-yx, \\
\lbrack x,y,z] & =x(yz)-(xy)z.
\end{align*}
These two multilinear maps provide an easy way for stating commutativity or
associativity of the algebra $A$: the algebra $A$ is commutative if and only
if the commutator is identically $0$, it is associative if the associator
map is $0$. The associator is also useful for characterizing alternativity.

\begin{proposition}
\label{th1/1}A $k$--algebra $A$ is alternative if and only if the associator
map $(x,y,z)\mapsto\lbrack x,y,z\rbrack$ is alternating.
\end{proposition}

In fact, the identities (\ref{eq1/20}) can be written equivalently
\begin{equation*}
\lbrack a,a,c]=0,%
\quad[c,a,a]=0.
\end{equation*}
They are obviously verified if the associator is alternating.

Conversely, let $A$ be alternative; then $[x,y,z]$ is $0$ for $x=y$ or $y=z$%
. This means that $[x,y,z]$ is alternating with respect to $(x,y)$ and with
respect to $(y,z) $. As the transpositions $(12)$ and $(23)$ generate the
symmetric group $\mathfrak{S}_{3}$, it follows that the associator is a
trilinear alternating map.

As a consequence, in an alternative algebra, we have
\begin{equation*}
\lbrack a,b,a]=0,
\end{equation*}
which can also be written
\begin{equation}
a(ba)=(ab)a.   \label{eq1/21}
\end{equation}
An algebra satisfying (\ref{eq1/21}) is called \index{flexible algebra}\emph{flexible.} In such an
algebra, we will simply write $aba$ for $a(ba)=(ab)a$. In a composition
algebra $(A,n)$, as $a+\widetilde{a}$ is a scalar (a multiple of $e$), the
identity (\ref{eq1/21}) is equivalent to
\begin{equation}
\widetilde{a}(ba)=(\widetilde{a}b)a.   \label{eq1/22}
\end{equation}

In an alternative algebra, we have the important
\index{MoufangR@\textsc{Moufang, Ruth} (1905--1977)}%
\index{Moufang identities}\emph{Moufang identities}:

\begin{theorem}[Ruth Moufang]
\label{th1/2}In an alternative algebra, the following identities are true:
\begin{align}
a(x(ay)) & =(axa)y,  \label{eq1/23} \\
((xa)ya) & =x(aya),  \label{eq1/24} \\
(ax)(ya) & =a(xy)a.   \label{eq1/25}
\end{align}
\end{theorem}

They are called respectively the \emph{left, right }and\emph{\ central}
Moufang identity.

\begin{proof}
From the definitions, we get
\begin{equation*}
a(x(ay))-(axa)y=[a,x,ay]+[ax,a,y];
\end{equation*}
the right hand side is symmetric in $(x,y)$, so it is enough to check that
it vanishes for $x=y$. We have $[a,x,ax]=[ax,a,x]$ by Proposition \ref{th1/1}%
; repeatedly using (\ref{eq1/20}), we obtain
\begin{equation*}
\lbrack a,ax,x]=a((ax)x)-(a(ax))x=a^{2}x^{2}-a^{2}x^{2}=0.
\end{equation*}
This proves the left identity (\ref{eq1/23}). The right identity (\ref%
{eq1/24}) is proved in the same way; we also note that it is just the left
identity in the opposite algebra, which is also alternative.

We also get from the definitions and from alternativity
\begin{align*}
(ax)(ya)-a(xy)a & =-[a,x,ya]+a[x,y,a] \\
& =[a,ya,x]+a[y,a,x] \\
& =a(y(ax))-((aya)x).
\end{align*}
The last expression vanishes by (\ref{eq1/23}), so the central identity (\ref%
{eq1/25}) is proved.
\end{proof}

The following proposition allows us to characterize composition algebras
among alternative algebras:

\begin{proposition}
\label{th1/3}Let $A$ be an algebra with unit element $e$. Assume there is an
involutive anti-automorphism $a\mapsto%
\widetilde{a}$ of $A$ (with $\widetilde{e}=e$) such that $a+\widetilde{a}$
and $a\widetilde{a}$ are scalars (multiples of $e$) for all $a\in A$. Define
$n:A\rightarrow k$ by
\begin{equation*}
n(a)=a\widetilde{a}.
\end{equation*}
Then $(A,n)$ is a composition algebra if and only if $A$ is alternative and $%
n$ is non-singular. In this case, the Cayley conjugation in $(A,n)$ is $%
a\mapsto\widetilde{a}$.
\end{proposition}

\begin{proof}
Let $a,b\in A$; then $a+\widetilde{a}=\alpha$, $b+\widetilde{b}=\beta$ with $%
\alpha,\beta\in ke$. We have then, using alternativity and the central
Moufang identity,
\begin{align*}
n(ab) & =(ab)\ (ab)\widetilde{\ }=(ab)(\widetilde{b}\widetilde{a}%
)=(ab)((\beta -b)(\alpha-a)) \\
& =ab\beta\alpha-ab^{2}\alpha-aba\beta+ab^{2}a \\
& =a(b(\beta-b))(\alpha-a)=a(b\widetilde{b})\widetilde{a}=n(b)a\widetilde {a}
\\
& =n(a)n(b).
\end{align*}
This shows that $n$ is multiplicative; if $n$ is non-singular, $\left(
A,n\right) $ is a composition algebra. The bilinear form associated to $n$
is then $\left( a:b\right) =a\widetilde{b}+b\widetilde{a}$, which shows that
the trace is $t(a)=\left( a:e\right) =a+\widetilde{a}$ and that $\widetilde{a%
}$ is indeed the Cayley conjugate of $a$ in $\left( A,n\right) $.
\end{proof}

\subsection{Cayley-Dickson extensions: analysis}

We are going to describe the
\index{CayleyA@\textsc{Cayley, Arthur} (1821--1895)}%
\index{DicksonL@\textsc{Dickson, Leonard Eugene} (1874--1954)}%
\index{Cayley--Dickson extension}\emph{Cay\-ley--Dickson extension process}:
start from the subalgebra $A_{0}=ke$. This process allows one to construct
successive subalgebras $A_{1},A_{2},A_{3}$, each time doubling the dimension
(as vector space) and terminates at most on the third step.

Let us first examine when a subalgebra $B$ of a composition algebra $\left(
A,n\right) $ is itself a composition algebra.

\begin{proposition}
\label{th1/4}Let $\left( A,n\right) $ and $\left( B,n^{\prime}\right) $ be
composition algebras with unit elements $e,e^{\prime}$. If $f:B\rightarrow A$
is an algebra homomorphism (with $fe^{\prime}=e$) and $f$ is injective, then
$f$ is a (partial) isometry:
\begin{equation*}
n(fx)=n^{\prime}(x)\quad(x\in B).
\end{equation*}
\end{proposition}

\begin{proof}
Let $x\in B$ and let $y=fx$. We have $x^{2}+t^{\prime}(x)x+n^{\prime
}(x)e^{\prime}=0$, which gives
\begin{equation*}
y^{2}+t^{\prime}(x)y+n^{\prime}(x)e=0;
\end{equation*}
comparing with $y^{2}+t(y)y+n(y)e=0$, we get
\begin{equation*}
\left( t(y)-t^{\prime}(x)\right) y+\left( n(y)-n^{\prime}(x)\right) e=0.
\end{equation*}
If $\left( y,e\right) $ is free, then $n(y)=n^{\prime}(x)$. If $y=\lambda e$%
, it follows from the injectivity of $f$ that $x=\lambda e^{\prime}$, and
then again $n(y)=n^{\prime}(x)=\lambda^{2}$.
\end{proof}

Proposition \ref{th1/4} shows that if $\left( B,n^{\prime}\right) $ is a
composition subalgebra of $\left( A,n\right) $, the norm of $B$ has to be
the restriction of the norm of $A$. If $\left( A,n\right) $ is a composition
algebra and $B$ is a subalgebra, then $B$ is a composition subalgebra of $%
\left( A,n\right) $ if and only if $\left. n\right| _{B}$ is non-singular,
that is, if $B\cap B^{\perp}=0$.

\begin{proposition}
\label{th1/5}Let $\left( A,n\right) $ be a composition algebra and let $B$ a
composition subalgebra. Assume $B\neq A$. Let $v\in B^{\perp}$ be
non-isotropic: $n(v)=-\mu\neq0$. Then

1)~the vector subspace $vB$ is orthogonal to $B$ and the map $\gamma
_{v}:x\mapsto vx$ is an isomorphism from $B$ onto $vB$;

2)~the subalgebra $C=\left\langle B,v\right\rangle $ generated by $B$ and $v$
is (as a vector space) $C=B\oplus vB$;

3)~$C$ is a composition subalgebra;

4)~the product in $C=B\oplus vB$ is defined by
\begin{equation}
\left( a_{1}+vb_{1}\right) \left( a_{2}+vb_{2}\right) =a_{1}a_{2}+\mu b_{2}%
\widetilde{b_{1}}+v\left( \widetilde{a_{1}}b_{2}+a_{2}b_{1}\right) ;
\label{eq1/26}
\end{equation}
the norm and the Cayley conjugation are defined by
\begin{align}
n(a+vb) & =n(a)-\mu n(b),  \label{eq1/27} \\
( a+vb)\widetilde{\ } & =\widetilde{a}-vb.   \label{eq1/28}
\end{align}
\end{proposition}

\begin{proof}
As $v\perp B$, we have in particular $\left( v:1\right) =0$, which implies $%
\widetilde{v}=-v$ and $v^{2}=-n(v)=\mu$. For each $b\in B$, we have
\begin{equation*}
0=\left( b:v\right) =\widetilde{b}v+\widetilde{v}b=\widetilde{b}v-vb,
\end{equation*}
which implies
\begin{equation}
vb =\widetilde{b}v,\quad ( vb)\widetilde{\ } =-vb.  \label{eq1/29}
\end{equation}
This proves (\ref{eq1/28}).

Let $a,b\in B$; then $( a:vb ) =( a\widetilde{b}:v) =0$. This shows $B\perp
vB$. As $B\cap B^{\perp}=0$, we have $B\cap vB=0$. The relation $v( vx) =\mu
x$ proves that $\gamma_{v}:x\mapsto vx$ is an isomorphism from $B$ onto $vB.$
The relation (\ref{eq1/27}) directly follows from $B\perp vB$ and $n(v)=-\mu$%
. It shows that the restriction of $n$ to $B\oplus vB$ is non-singular. It
remains to prove the relations (\ref{eq1/26}), which will imply that $%
C=\left\langle B,v\right\rangle =B\oplus vB$ and that $C$ is a composition
subalgebra.

Using (\ref{eq1/29}), the central Moufang identity and alternativity, we get
\begin{equation*}
( vb_{1}) ( vb_{2}) =( vb_{1}) ( \widetilde{b_{2}}v) =v( b_{1}\widetilde{%
b_{2}}) v=v^{2}b_{2}\widetilde{b_{1}}=\mu b_{2}\widetilde{b_{1}}.
\end{equation*}
Using the left Moufang identity and (\ref{eq1/29}), we have
\begin{equation*}
v\left( a_{1}\left( vb_{2}\right) \right) =\left( va_{1}v\right) b=\left(
v^{2}\widetilde{a_{1}}\right) b_{2}=\mu\widetilde{a_{1}}b_{2}
\end{equation*}
and, multiplying by $\mu^{-1}v$,
\begin{equation*}
a_{1}\left( vb_{2}\right) =v\left( \widetilde{a_{1}}b_{2}\right) .
\end{equation*}
Conjugating (after $b_{2}\leftarrow b_{1}$ and $a_{1}\leftarrow\widetilde {%
a_{2}}$), we obtain $( \widetilde{b_{1}}v) a_{2}=( \widetilde{b_{1}}%
\widetilde{a_{2}}) v$ and, using again (\ref{eq1/29}),
\begin{equation*}
\left( vb_{1}\right) a_{2}=v\left( a_{2}b_{1}\right) .
\end{equation*}
\end{proof}

\subsection{Cayley-Dickson extensions: construction}

Let $\left( A,n\right) $ be a composition algebra. Let $A^{\prime}=A\times A
$. In view of Proposition \ref{th1/5}, we consider on the vector space $%
A^{\prime}$ the product defined by
\begin{equation*}
( a_{1},b_{1}) ( a_{2},b_{2}) =( a_{1}a_{2}+\mu b_{2}\widetilde{b_{1}},%
\widetilde{a_{1}}b_{2}+a_{2}b_{1})
\end{equation*}
and the quadratic form $n^{\prime}$ defined by
\begin{equation*}
n^{\prime}(a,b)=n(a)-\mu n(b);
\end{equation*}
we ask whether $\left( A^{\prime},n^{\prime}\right) $ is a composition
algebra. In this case, it also follows from Proposition \ref{th1/5} that the
conjugation in $\left( A^{\prime},n^{\prime}\right) $ will be given by
\begin{equation*}
( a,b)\widetilde{\ }=( \widetilde{a},-b) .
\end{equation*}
With these definitions, $A$ can be identified with the subalgebra $A\times0$
of $A^{\prime}$, by $a\mapsto\left( a,0\right) $ ; the norm and conjugation
in $A^{\prime}$ then extend those of $A$. If we set $v=\left( 0,e\right) ,$
then $v\left( b,0\right) =\left( 0,b\right) .$ So we can write $A^{\prime
}=A\oplus vA$ and the operation rules in $A^{\prime}$ are given by (\ref%
{eq1/26}), (\ref{eq1/27}), (\ref{eq1/28}). It is easily seen that the
conjugation in $A^{\prime}$ is an involutive antiautomorphism. Moreover, we
have $\left( a,b\right) +( a,b)\widetilde{\ }=t(a)\left( e,0\right) $ and $%
\left( a,b\right) ( a,b)\widetilde{\ }=\left( a,b\right) \left( \widetilde{a}%
,-b\right) =n^{\prime}(a,b)\left( e,0\right) $. Also, the definition of $%
n^{\prime}$ shows that it is non-singular if $\mu\neq0$. So all conditions
of Proposition \ref{th1/3} are fulfilled by $A^{\prime}$ with the product (%
\ref{eq1/26}) and the conjugation (\ref{eq1/28}), except alternativity. The
answer to this last question is given by the following proposition.

\begin{proposition}
\label{th1/6}Let $\left( A,n\right) $ be a composition algebra. For $\mu
\neq0$, denote by $A(\mu)$ the algebra
\begin{equation*}
A(\mu)=A\oplus vA
\end{equation*}
with the product
\begin{equation*}
\left( a_{1}+vb_{1}\right) \left( a_{2}+vb_{2}\right) =a_{1}a_{2}+\mu b_{2}%
\widetilde{b_{1}}+v\left( \widetilde{a_{1}}b_{2}+a_{2}b_{1}\right) .
\end{equation*}
Then

1)~$A(\mu)$ is commutative if and only if $A=ke$;

2)~$A(\mu)$ is associative if and only if $A$ is associative and commutative;

3)~$A(\mu)$ is alternative if and only if $A$ is associative.
\end{proposition}

\begin{proof}
1)~The definition of the product in $A(\mu)$ implies $av=v\widetilde{a}$ for
all $a\in A$. If $A(\mu)$ is commutative, we have $va=v\widetilde{a}$, which
implies $a=\widetilde{a}$ and $a\in ke$ for all $a\in A$. This shows $A=ke$.

Conversely, the algebra $k(\mu)$ with the product
\begin{equation*}
\left( a_{1}+vb_{1}\right) \left( a_{2}+vb_{2}\right) =a_{1}a_{2}+\mu
b_{1}b_{2}+v\left( a_{1}b_{2}+a_{2}b_{1}\right)
\end{equation*}
is clearly commutative (and associative).

2)~If $A(\mu)$ is associative, $A$ is also associative. For $a,b\in A$, it
is easily checked that
\begin{equation*}
\lbrack a,b,v\rbrack =v\lbrack \widetilde{a},\widetilde{b}\rbrack ;
\end{equation*}
so $A(\mu)$ associative implies $A$ commutative.

Conversely, assume $A$ is associative; then routine computations using the
definitions show that, for $x,y,z\in A$, one has
\begin{align*}
&\left[ vx,y,z\right] =v\left( \left[ y,z\right] x\right) ,\quad\left[ x,vy,z%
\right] =v\left( \left[ \widetilde{x},z\right] y\right) ,\quad\left[ x,y,vz%
\right] =v\left( \left[ \widetilde {x},\widetilde{y}\right] z\right) , \\
&\left[ x,vy,vz\right] =\mu\left[ x,z\widetilde{y}\right] ,\quad\left[
vx,y,vz\right] =\mu\left[ \widetilde{y},z\widetilde{x}\right] ,\quad\left[
vx,vy,z\right] =\mu\left[ z,y\widetilde{x}\right] , \\
&\left[ vx,vy,vz\right] =\mu v\left( \left[ z\widetilde{y},x\right] +x\left[
z,\widetilde{y}\right] \right) .
\end{align*}
If moreover $A$ is commutative, we see that $A(\mu)$ is associative.

3)~Let $x,y\in A(\mu)$. As $x+\widetilde{x}$ is a scalar, we have
\begin{equation*}
\left[ x,x,y\right] =-\left[ x,\widetilde{x},y\right]
\end{equation*}
and
\begin{equation*}
\left[ y,x,x\right] \widetilde{\ }=-\left[ y,x,\widetilde{x}\right]
\widetilde{\ }=\left[ x,\widetilde{x},y\right] .
\end{equation*}
This shows that $A(\mu)$ is alternative if and only if $\left[ x,\widetilde {%
x},y\right] =0$ for all $x,y\in A(\mu)$. Now assume $A$ is alternative and
let $x=x_{1}+vx_{2}$, $y=y_{1}+y_{2}$; then
\begin{equation*}
\left[ x,\widetilde{x},y\right] =-\mu\left[ x_{1},y_{2},\widetilde{x_{2}}%
\right] +v\left[ \widetilde{x_{1}},y_{1},x_{2}\right] ,
\end{equation*}
which shows that $A(\mu)$ is alternative if and only if $A\,$is associative.
\end{proof}

\begin{theorem}
\label{th1/7}A composition algebra is (as a vector space) of dimension $1$, $%
2$, $4$ or $8$.
\end{theorem}

\begin{proof}
Let $A$ be a composition algebra over $k$. Let $A_{0}=ke$. If $A\neq A_{0}$,
there exists $v_{1}\perp e$, with $n(v_{1})=-\mu_{1}\neq0$ ; the composition
subalgebra $A_{1}=ke\oplus kv_{!}$ is commutative and associative. If $A\neq
A_{1}$, there exists $v_{2}\perp A_{1}$, with $n(v_{2})=-\mu_{2}\neq0$ ; the
composition subalgebra $A_{2}=A_{1}\oplus v_{2}A_{1}$ is associative, but
not commutative, of dimension $4$. If $A\neq A_{2},$ there exists $%
v_{3}\perp A_{2}$, with $n(v_{3})=-\mu_{3}\neq0$ ; the composition
subalgebra $A_{3}=A_{2}\oplus v_{3}A_{2}$ is alternative (as $A_{2}$ is
associative), but not associative (as $A_{2}$ is not commutative), of
dimension $8$. Then $A=A_{3}$, as $A_{3}$ is not associative.
\end{proof}

\subsection{Classification of composition algebras over $\mathbb{R}$ or $%
\mathbb{C}$}

We consider the composition algebras
\begin{equation*}
k,\quad k(\mu_{1}),\quad k(\mu_{1},\mu_{2})=\left( k(\mu_{1})\right)
(\mu_{2}),\quad k(\mu_{1},\mu_{2},\mu_{3})=\left( k(\mu_{1},\mu_{2})\right)
(\mu_{3}),
\end{equation*}
for non-zero $\mu_{1},\mu_{2},\mu_{3}\in k$. It follows from the proof of
Theorem \ref{th1/7} that each composition algebra is isomorphic to one of
these for a suitable choice of $\mu_{1},\mu_{2},\mu_{3}$. We want to make
this statement more precise. First, we show that, if the norms of two
composition algebras are linearly equivalent, these composition algebras are
isomorphic.

\begin{proposition}
\label{th1/8}Let $\left( A,n\right) $ and $\left( A^{\prime},n^{\prime
}\right) $ be composition algebras. Then $A$ and $A^{\prime}$ are isomorphic
(as unital algebras) if and only there exists a linear isomorphism $%
f:A\rightarrow A^{\prime}$ such that $n^{\prime}\circ f=n$.
\end{proposition}

\begin{proof}
By Proposition \ref{th1/4}, an isomorphism of composition algebras preserves
norms.

Assume there exists a linear isomorphism $f:A\rightarrow A^{\prime}$ such
that $n^{\prime}\circ f=n$. Let $B,B^{\prime}$ be proper composition
subalgebras of $A,A^{\prime}$ respectively, such that there exists an
algebra isomorphism $g:B\rightarrow B^{\prime}$; then $n^{\prime}\circ
g=\left. n\right| _{B}$ by Proposition \ref{th1/4}. By Witt's theorem, $g$
can be extended to a vector space isomorphism $h:A\rightarrow A^{\prime}$
such that $n^{\prime}\circ h=n$. Let $v\in B^{\perp}$ such that $%
n(v)=-\mu\neq0$; take $v^{\prime}=h(v)$, which implies $v^{\prime}\in
B^{\prime\perp}$ and $n^{\prime}(v^{\prime})=n(v)=-\mu $. Let
\begin{equation*}
\widehat{g}:B\oplus vB\rightarrow B^{\prime}\oplus v^{\prime}B^{\prime}
\end{equation*}
be defined by
\begin{equation*}
\widehat{g}\left( a+vb\right) =g(a)+v^{\prime}g(b).
\end{equation*}
By Proposition \ref{th1/5}, $\widehat{g}$ is then an algebra isomorphism
between the composition subalgebras $B\oplus vB$ and $B^{\prime}\oplus
v^{\prime}B^{\prime}$. Starting from the trivial isomorphism $%
g_{0}:ke\rightarrow ke^{\prime}$ and iterating this process at most three
times, we get an algebra isomorphism from $A$ onto $A^{\prime}$.
\end{proof}

Assume that the ground field is $k=\mathbb{C}$. In each dimension, there is
only one non-singular quadratic form, up to linear equivalence. So
Proposition \ref{th1/8} implies:

\begin{theorem}
\label{th1/9}On $k=\mathbb{C}$, there exist, up to isomorphism, exactly four
composition algebras $A_{j}$ $(0\leq j\leq3)$, of respective dimension $2^{j}
$.
\end{theorem}

Assume now that the ground field is $k=\mathbb{R}$. In this case,
non-singular quadratic forms are classified, up to linear equivalence, by
their signature. The signature for $A_{0}^{+}=\mathbb{R}$ is $\left(
1,0\right) $. Let us show that for other composition algebras over $\mathbb{R%
}$, the signature needs to be $\left( 2a,0\right) $ or $\left( a,a\right) $.
If $\left( B,n\right) $ is a composition algebra, we know that the norm of
the Cayley-Dickson extension $B(\mu)$ is given by $n^{\prime}(a+vb)=n(a)-\mu
n(b)$. The signature of $n^{\prime}$ is

\begin{itemize}
\item $\left( 2a,0\right) $ if the signature of $n$ is $\left( a,0\right) $
and $\mu<0$;

\item $\left( a,a\right) $ if the signature of $n$ is $\left( a,0\right) $
and $\mu>0$;

\item $\left( 2b,2b\right) $ if the signature of $n$ is $\left( b,b\right) $.
\end{itemize}

The assumption on the signature can then be proved by induction.

\begin{theorem}
\label{th1/10}On $k=\mathbb{R}$, there exist, up to isomorphism, seven
composition algebras:

the ``compact'' algebras $A_{j}^{+}$ $(0\leq j\leq3)$ of dimension $2^{j}$,
with positive-definite norm;

the ``split'' algebras $A_{j}^{-}$ $(1\leq j\leq3)$ of dimension $2^{j}$ and
signature $\left( 2^{j-1},2^{j-1}\right) $.
\end{theorem}

When $k=\mathbb{C}$, a model for the composition algebra of dimension $4$ is
$A_{2}\cong\mathcal{M}_{2,2}(\mathbb{C})$ ($2\times2$ complex matrices),
with the determinant as norm; a model for $A_{1}$ is the subalgebra of
diagonal $2\times2$ complex matrices. The non-associative composition
algebra $A_{3}$ is called the
\index{Cayley algebra!complex ---}\emph{complex Cayley algebra} or the
\index{octonions!complex ---}\emph{algebra of
complex octonions}. It can be constructed, for example, as $A_{2}\left(
-1\right) $; but this is in most cases irrelevant and it will be more
important to know that this composition algebra of dimension $8$ exists and
is unique up to isomorphism. The algebra $A_{3}$ will be denoted by
\index{$OC@$\mathbb{O}_{\mathbb{C}}$}$\mathbb{O}_{\mathbb{C}}$.

In the case $k=\mathbb{R}$, models for the compact composition algebras of
dimension $2$ and $4$ are respectively $A_{1}^{+}\cong\mathbb{C}$ (with norm
$n(z)=\left| z\right| ^{2}$) and $A_{2}^{+}=\mathbb{H}$ \index{$H@$\mathbb{H}$}(the field of
quaternions), which can be described as
\begin{equation*}
\mathbb{H}=\left\{ q=
\begin{pmatrix}
a & -%
\overline{b} \\
b & \overline{a}%
\end{pmatrix}
;\ a,b\in\mathbb{C}\right\} ,
\end{equation*}
with norm $n(q)=a\overline{a}+b\overline{b}$. The compact non-associative
real composition algebra $A_{3}^{+}$ is known as the \emph{algebra of Cayley
numbers}, the \index{octonions}\emph{algebra of octonions} or the
\index{Cayley algebra!compact ---}\emph{Cayley real division algebra}.
It will be denoted by \index{$O@$\mathbb{O}$, $\mathbb{O}_{c}$}$\mathbb{O}$ or
$\mathbb{O}_{c}$; it can be constructed as $\mathbb{H}\left(
-1\right) $. Again the most important point is that $\mathbb{O}$ is a real
composition algebra of dimension $8$ with positive norm, and is unique up to
isomorphism.

The split composition algebras $A_{1}^{-}$ and $A_{2}^{-}$ are respectively
isomorphic to the algebra of diagonal $2\times2$ real matrices and to the
algebra $\mathcal{M}_{2,2}(\mathbb{R})$ of $2\times2$ real matrices, with
the determinant as norm. The algebra $A_{3}^{-}$ can be constructed as $%
\mathbb{R}\left( 1,1,1\right) $; the signature of its norm is $\left(
4,4\right) $. It is denoted by \index{$Os@$\mathbb{O}_{s}$}%
$\mathbb{O}_{s}$ and called the \index{Cayley algebra!split ---}\emph{split Cayley algebra}.

The real composition algebras can be \emph{complexified} in a natural way.
The complexification is then isomorphic to the complex composition algebra
of the corresponding dimension.


\section{Exceptional Jordan triple systems}



\subsection{The space $H_{3}(\mathbb{O})$}

In this section, $\mathbb{O}$ denotes a Cayley algebra over $k=\mathbb{R}$
or $\mathbb{C}$.

\begin{definition}
\label{def2/1}The space \index{$H30@$\mathcal{H}_{3}(\mathbb{O}_{\mathbb{C}})$}%
$H_{3}(\mathbb{O})$ is the $k$-vector space (with
the natural operations) of $3\times3$ matrices with entries in $\mathbb{O}$,
which are Hermitian with respect to the Cayley conjugation in $\mathbb{O}$.
\end{definition}

An element $a\in H_{3}(\mathbb{O})$ will be written
\begin{equation}
a=
\begin{pmatrix}
\alpha_{1} & a_{3} & \widetilde{a_{2}} \\
\widetilde{a_{3}} & \alpha_{2} & a_{1} \\
a_{2} & \widetilde{a_{1}} & \alpha_{3}%
\end{pmatrix},   \label{eq2/1}
\end{equation}
with $\alpha_{1},\alpha_{2},\alpha_{3}\in k$ and $a_{1},a_{2},a_{3}\in%
\mathbb{O}$. Instead of (\ref{eq2/1}), we will also write
\begin{equation}
a=\sum_{j=1}^{3}\alpha_{j}e_{j}+\sum_{j=1}^{3}F_{j}(a_{j}),   \label{eq2/50}
\end{equation}
with the obvious definitions for $e_{j}$ and $F_{j}(a_{j})$. The vector
space $H_{3}(\mathbb{O})$ decomposes into the direct sum
\begin{equation}
H_{3}(\mathbb{O})=ke_{1}\oplus ke_{2}\oplus ke_{3}\oplus\mathcal{F}_{1}\oplus%
\mathcal{F}_{2}\oplus\mathcal{F}_{3},   \label{eq2/2}
\end{equation}
where $\mathcal{F}_{j}=\left\{ F_{j}(a)\mid a\in\mathbb{O}\right\} $. The
subspaces $\mathcal{F}_{j}$ are $8$-dimensional and
\begin{equation*}
\dim_{k}H_{3}(\mathbb{O})=27.
\end{equation*}

On $H_{3}(\mathbb{O})$, define a bilinear form by
\begin{equation}
(a:b)=\sum_{j=1}^{3}\alpha_{j}\beta_{j}+\sum_{j=1}^{3}(a_{j}:b_{j})
\label{eq2/3}
\end{equation}
for
\begin{align*}
a & =\sum_{j=1}^{3}\alpha_{j}e_{j}+\sum_{j=1}^{3}F_{j}(a_{j}), \\
b & =\sum_{j=1}^{3}\beta_{j}e_{j}+\sum_{j=1}^{3}F_{j}(b_{j});
\end{align*}
in (\ref{eq2/3}), $(a_{j}:b_{j})$ denotes the scalar product in $\mathbb{O}$%
. The form defined by (\ref{eq2/3}) is clearly non-singular and the
decomposition (\ref{eq2/2}) is orthogonal with respect to it. We will refer
to $(a:b)$ as the \emph{scalar product} of $a$ and $b$ in $H_{3}(\mathbb{O})$%
.

\begin{definition}
\label{def2/2}The \index{$a@$a^{\#}$, adjoint of $a\in H_{3}(\mathbb{O})$}%
\emph{adjoint} $a^{\#}$ of an element $a\in H_{3}(\mathbb{O%
})$, written in the form (\ref{eq2/50}), is defined by
\begin{equation}
a^{\#}=\sum_{i}\left( \alpha_{j}\alpha_{k}-n(a_{i})\right) e_{i}+\sum _{i}%
\widetilde{F_{i}}\left( a_{j}a_{k}-\alpha_{i}\widetilde{a_{i}}\right) .
\label{eq2/4}
\end{equation}
\end{definition}

In (\ref{eq2/4}) and below, $\sum_{i}$ means $\sum_{i=1}^{3}$ and $j,k$ are
defined by $(i,j,k)$ being an \emph{even} permutation of $(1,2,3)$; $%
\widetilde{F_{i}}(c)$ stands for $F_{i}(\widetilde{c})$.

\begin{definition}
\label{def2/3}The symmetric bilinear map, associated to the quadratic map $%
a\mapsto a^{\#}$, is called the
\index{FreudenthalH@\textsc{Freudenthal, Hans} (1905--1990)}%
\index{Freudenthal product}\emph{Freudenthal product}. The Freudenthal
product of $a,b\in H_{3}(\mathbb{O})$ is denoted $a \times b$ and is defined
by
\begin{equation*}
a\times b=(a+b)^{\#}-a^{\#}-b^{\#},\quad a\times a=2a^{\#},
\end{equation*}
\end{definition}

It follows directly from the definitions that%
\begin{align}
a\times b= & \sum_{i}\left(
\alpha_{j}\beta_{k}+\alpha_{k}\beta_{j}-(a_{i}:b_{i})\right) e_{i}  \notag \\
& +\sum_{i}%
\widetilde{F_{i}}( a_{j}b_{k}+b_{j}a_{k}-\alpha _{i}\widetilde{b_{i}}%
-\beta_{i}\widetilde{a_{i}}) .
\end{align}
The following multiplication rules hold:%
\begin{align}
e_{i}\times e_{i} & =0,\quad e_{i}\times e_{j}=e_{k},  \notag \\
e_{i}\times F_{i}(b) & =-F_{i}(b),\quad e_{i}\times F_{j}(b)=0,\quad
\label{eq2/27} \\
F_{i}(a)\times F_{i}(b) & =-(a:b)e_{i},\quad F_{i}(a)\times F_{j}(b)=%
\widetilde{F_{k}}(ab),  \notag
\end{align}
where $(i,j,k)$ always stands for an even permutation of $(1,2,3)$.

\begin{proposition}
\label{th2/1}
\begin{equation}
\left( a\times b:c\right) =\left( a:b\times c\right) \qquad(a,b,c\in H_{3}(%
\mathbb{O})).   \label{eq2/5}
\end{equation}
\end{proposition}

\begin{proof}
For
\begin{align*}
a & =\sum_{j=1}^{3}\alpha_{j}e_{j}+\sum_{j=1}^{3}F_{j}(a_{j}), \quad
b =\sum_{j=1}^{3}\beta_{j}e_{j}+\sum_{j=1}^{3}F_{j}(b_{j}), \\
c & =\sum_{j=1}^{3}\gamma_{j}e_{j}+\sum_{j=1}^{3}F_{j}(c_{j}),
\end{align*}
by applying the definitions we obtain
\begin{align*}
\left( a\times b:c\right)= &\sum_{i}\left(
\alpha_{j}\beta_{k}+\alpha_{k}\beta_{j}-(a_{i}:b_{i})\right) \gamma_{i} \\
& +\sum_{i}\left( a_{j}b_{k}+b_{j}a_{k}-\alpha_{i}\widetilde{b_{i}}-\beta_{i}%
\widetilde{a_{i}}:\widetilde{c_{i}}\right) \\
=&\sum_{(i,j,k)\in\mathfrak{S}_{3}}\alpha_{i}\beta_{j}\gamma_{k}+%
\sum_{(i,j,k)\in\mathfrak{S}_{3}}t(a_{i},b_{j},c_{k}) \\
& -\sum_{i}\left(
(a_{i}:b_{i})\gamma_{i}+(b_{i}:c_{i})\alpha_{i}+(c_{i}:a_{i})\beta_{i}%
\right) .
\end{align*}
(Recall that $t(x,y,z)=t((xy)z)=(xy:\widetilde{z})$ for $x,y,z\in\mathbb{O}$%
). This shows that $\left( a\times b:c\right) $ is symmetric with respect to
$(a,b,c)$.
\end{proof}

\begin{definition}
\label{def2/4}Let $T$ denote the trilinear symmetric form on $H_{3}(\mathbb{O%
})$ defined by
\begin{equation*}
T(a,b,c)=\left( a\times b:c\right) .
\end{equation*}
The
\index{determinant in $H_{3}(\mathbb{O})$}\emph{determinant} in $H_{3}(%
\mathbb{O})$ is the associated polynomial of degree $3$, defined by
\begin{equation}
\det a=%
\frac{1}{3!}T(a,a,a)=\frac{1}{3}(a^{\#}:a).   \label{eq2/6}
\end{equation}
\end{definition}

From the expression of $\left( a\times b:c\right) $, we deduce
\begin{equation}
\det
a=\alpha_{1}\alpha_{2}\alpha_{3}-\sum_{i}%
\alpha_{i}n(a_{i})+a_{1}(a_{2}a_{3})+(\widetilde{a_{3}}\widetilde{a_{2}})%
\widetilde{a_{1}}.   \label{eq2/7}
\end{equation}
This relation may also be taken as a definition of $\det a$. It is an
extension of the classical ``Sarrus' rule" for $3\times3$ matrices, but with
suitable parentheses in products like $a_{1}(a_{2}a_{3})$, due to the
non-associativity of the Cayley algebra.

\begin{proposition}
\label{th2/2}
\begin{equation}
(a^{\#})^{\#}=(\det a)a.   \label{eq2/8}
\end{equation}
\end{proposition}

\begin{proof}
Let
\begin{align*}
& a  =\sum_{j=1}^{3}\alpha_{j}e_{j}+\sum_{j=1}^{3}F_{j}(a_{j}), \quad
a^{\#} =\sum_{j=1}^{3}\beta_{j}e_{j}+\sum_{j=1}^{3}F_{j}(b_{j}), \\
& (a^{\#})^{\#}  =\sum_{j=1}^{3}\gamma_{j}e_{j}+\sum_{j=1}^{3}F_{j}(c_{j}).
\end{align*}
From the definition (\ref{eq2/4}) and the properties of Cayley algebras, we
get%
\begin{align*}
\gamma_{i} & =\beta_{j}\beta_{k}-n(b_{i}) \\
& =\left( \alpha_{k}\alpha_{i}-n(a_{j})\right) \left( \alpha_{i}\alpha
_{j}-n(a_{k})\right) -n\left( a_{j}a_{k}-\alpha_{i}\widetilde{a_{i}}\right)
\\
& =\alpha_{i}^{2}\alpha_{j}\alpha_{k}-\alpha_{i}\alpha_{j}n(a_{j})-\alpha
_{i}\alpha_{k}n(a_{k})+n(a_{j})n(a_{k}) \\
& \quad -n\left( a_{j}a_{k}\right) -\alpha_{i}^{2}n(a_{j})+\alpha_{i}\left(
a_{j}a_{k}:\widetilde{a_{i}}\right) \\
& =\alpha_{i}\det a
\end{align*}
(using namely $n(a_{j})n(a_{k})=n\left( a_{j}a_{k}\right) $) and
\begin{align*}
c_{i} & =\widetilde{b_{k}}\widetilde{b_{j}}-\beta_{i}b_{i} \\
& =\left( a_{i}a_{j}-\alpha_{k}\widetilde{a_{k}}\right) \left(
a_{k}a_{i}-\alpha_{j}\widetilde{a_{j}}\right) -\left( \alpha_{j}\alpha
_{k}-n(a_{i})\right) \left( \widetilde{a_{k}}\widetilde{a_{j}}-\alpha
_{i}a_{i}\right) \\
&
=(a_{i}a_{j})(a_{k}a_{i})-\alpha_{k}n(a_{k})a_{i}-%
\alpha_{j}n(a_{j})a_{i}+n(a_{i})\widetilde{a_{k}}\widetilde{a_{j}} \\
&\quad +\alpha_{i}\alpha_{j}\alpha_{k}a_{i}-\alpha_{i}n(a_{i})a_{i} \\
& =\left( \det a\right) a_{i}
\end{align*}
(here we used the central Moufang identity
\begin{equation*}
(a_{i}a_{j})(a_{k}a_{i})=\left( a_{i}(a_{j}a_{k})\right) a_{i}
\end{equation*}
and $n(a_{i})\widetilde{a_{k}}\widetilde{a_{j}}=\left( (\widetilde{a_{k}}%
\widetilde{a_{j}})\widetilde{a_{i}}\right) a_{i}$).
\end{proof}

\begin{proposition}
\label{th2/3}The following identities hold in $H_{3}(\mathbb{O})$:
\begin{align}
\det(a^{\#}) & =(\det a)^{2};  \label{eq2/9} \\
\mathrm{d}(\det a).b & =\left( a^{\#}:b\right) ;  \label{eq2/10} \\
a^{\#}\times(a\times b) & =(\det a)b+\left( a^{\#}:b\right) a;
\label{eq2/11} \\
\left( a\times b:a^{\#}\times c\right) & =(\det a)(b:c)+\left(
a^{\#}:b\right) (a:c);  \label{eq2/12} \\
a\times\left( a^{\#}\times c\right) & =(\det a)c+(a:c)a^{\#};
\label{eq2/13} \\
(a\times b)\times(a\times c) & +a^{\#}\times(b\times c) = \notag \\
& =\left( a^{\#}:b\right) c+\left( a^{\#}:c\right) b+T(a,b,c)a;
\label{eq2/14} \\
a\times\left( \left( a\times b\right) \times c\right) & +b\times\left(
a^{\#}\times c\right)=  \notag \\
& =\left( a^{\#}:b\right) c+(b:c)a^{\#}+(a:c)a\times b;   \label{eq2/14a} \\
a^{\#}\times b^{\#}+(a\times b)^{\#} & =\left( a^{\#}:b\right) b+\left(
b^{\#}:a\right) a;  \label{eq2/15} \\
\left( a\times b^{\#}:a^{\#}\times b\right) & =3\det a\det b+(a:b)\left(
a^{\#}:b^{\#}\right) .   \label{eq2/16}
\end{align}
\end{proposition}

\begin{proof}
We have
\begin{equation*}
\det(a^{\#})=\frac{1}{3}\left( a^{\#}:(a^{\#})^{\#}\right) =\frac{1}{3}%
(a^{\#}:\det a\ a)=\left( \det a\right) ^{2}.
\end{equation*}
By differentiating the relation $\det a=\frac{1}{6}T(a,a,a)$, we get
\begin{equation*}
\mathrm{d}(\det a).b=\frac{1}{2}T(a,a,b)=\left( a^{\#}:b\right) .
\end{equation*}
We obtain (\ref{eq2/11}) and (\ref{eq2/14}) by successive differentiations
of (\ref{eq2/8}). The identity (\ref{eq2/12}) is obtained from (\ref{eq2/11}%
) by taking the scalar product with $c$ and using (\ref{eq2/5}). Using (\ref%
{eq2/5}) again and the fact that $\left( \ :\ \right) $ is non-singular, we
deduce (\ref{eq2/13}) from (\ref{eq2/12}). The relation (\ref{eq2/14a}) is
obtained by differentiating (\ref{eq2/13}). The relation (\ref{eq2/15}) is (%
\ref{eq2/14}) with $b=c$, and the identity (\ref{eq2/16}) is (\ref{eq2/12})
with $b=c^{\#}$.
\end{proof}

\subsection{The Hermitian Jordan triple system $H_{3}(\mathbb{O})$}

For the definition and general properties of
\index{Jordan triple system}Jordan triple systems, we refer the reader to
\cite{Loos1975}, \cite{Loos1977}, \cite{Roos1999}.

Let $\mathbb{O}_{c}$ be the compact Cayley algebra over $\mathbb{R}$, with
norm $n$ and Euclidean associated scalar product $\left( \ :\ \right) $. We
consider the complex Cayley algebra $\mathbb{O}$ as the complexification of $%
\mathbb{O}_{c}$: $\mathbb{O=C\otimes}_{\mathbb{R}}\mathbb{O}_{c};$ the
product, the Cayley conjugation and the norm on $\mathbb{O}$ are defined as
the natural extensions of those on $\mathbb{O}_{c}$: $(\alpha\otimes
a)(\beta\otimes b)=\alpha\beta\otimes ab$, $%
\widetilde{(\alpha\otimes a)}=\alpha\otimes\widetilde{a}$ and $%
n(\alpha\otimes a)=\alpha^{2}n(a)$ ($\alpha,\beta\in\mathbb{C}$, $a,b\in%
\mathbb{O}_{c}$). In addition, the algebra $\mathbb{O}$ has a complex
conjugation with respect to its ``real form'' $\mathbb{O}_{c}$, defined by
\begin{equation*}
\overline{(\alpha\otimes a)}=\overline{\alpha}\otimes a\qquad(\alpha \in%
\mathbb{C},\ a\in\mathbb{O}_{c});
\end{equation*}
this conjugation is antilinear and satisfies $\overline{ab}=\overline {a}%
\overline{b}$, in contrast with the Cayley conjugation which is complex
linear and satisfies $\widetilde{ab}=\widetilde{b}\widetilde{a}$.

The space $H_{3}(\mathbb{O}),$ with the operations defined in the previous
section, is then the complexification of the space $H_{3}(\mathbb{O}_{c})$
with the same operations. If
\begin{equation*}
a=\sum_{j=1}^{3}\alpha_{j}e_{j}+\sum_{j=1}^{3}F_{j}(a_{j})\in H_{3}(\mathbb{O%
}),
\end{equation*}
its complex conjugate with respect to $H_{3}(\mathbb{O}_{c})$ is defined by
\begin{equation*}
\overline{a}=\sum_{j=1}^{3}\overline{\alpha_{j}}e_{j}+\sum_{j=1}^{3}F_{j}(%
\overline{a_{j}}).
\end{equation*}
Clearly we have%
\begin{equation*}
\overline{a}^{\#}=\overline{a^{\#}},\quad\overline{a}\times\overline {b}=%
\overline{a\times b},\quad\det\overline{a}=\overline{\det a}.
\end{equation*}
On $\mathbb{O}$ and $H_{3}(\mathbb{O)},$ we define the Hermitian scalar
product
\begin{equation*}
\left( a\mid b\right) =( a:\overline{b}) .
\end{equation*}

\begin{definition}
\label{def2/5}The triple product $\{xyz\}$ on $H_{3}(\mathbb{O)}$, and the
related operators $Q$ and $D$, are defined by
\begin{align}
Q(x)y & =(x\mid y)x-x^{\#}\times\overline{y},  \label{eq2/17} \\
D(x,y)z & =\{xyz\}=(x\mid y)z+(z\mid y)x-(x\times z)\times\overline{y}.
\label{eq2/18}
\end{align}
\end{definition}

\begin{proposition}
\label{th2/4}With this triple product, $H_{3}(\mathbb{O})$ is a Hermitian
Jordan triple system.
\end{proposition}

\begin{proof}
The triple product defined by (\ref{eq2/18}) is clearly $\mathbb{C}$%
-bilinear symmetric in $(x,z)$ and antilinear in $y$. We are going to prove
that it satisfies the defining identities (J1) and (J2) of a Jordan triple
system.

Let us prove
\begin{equation}
D(x,y)Q(x)=Q(x)D(y,x).   \tag{J1}
\end{equation}
We have\pagebreak[0]%
\begin{align*}
D(x,y)Q(x)u & =(x\mid y)Q(x)u+(Q(x)u\mid y)x-(x\times Q(x)u)\times \overline{%
y} \\
& =(x\mid y)((x\mid u)x-x^{\#}\times\overline{u})
 +((x\mid u)x-x^{\#}\times\overline{u}\mid y)x \\
& \quad -(x\times((x\mid u)x-x^{\#}\times\overline{u}))\times\overline{y} \\
& =2(x\mid y)(x\mid u)x-(x\mid y)x^{\#}\times\overline{u}-(x^{\#}\mid
y\times u)x \\*
& \quad -2(x\mid u)x^{\#}\times\overline{y}+(x\times(x^{\#}\times\overline {u}%
))\times\overline{y};
\end{align*}
using (\ref{eq2/13}): $x\times\left( x^{\#}\times\overline{u}\right) =(\det
x)\overline{u}+(x\mid u)x^{\#}$, we get
\begin{align*}
D(x,y)Q(x)u & =2(x\mid y)(x\mid u)x-(x\mid y)x^{\#}\times\overline{u}-(x\mid
u)x^{\#}\times\overline{y} \\
& \quad -(x^{\#}\mid y\times u)x+(\det x)\overline{y}\times\overline{u}.
\end{align*}
In the same way,
\begin{align*}
Q(x)D(y,x)u & =(x\mid D(y,x)u)x-x^{\#}\times\overline{D(y,x)u} \\
& =(x\mid(y\mid x)u+(u\mid x)y-(y\times u)\times\overline{x})x \\
& \quad -x^{\#}\times\left( (x\mid y)\overline{u}+(x\mid u)\overline{y}-\overline{%
(y\times u)}\times x\right) \\
& =2(x\mid y)(x\mid u)x-2(x^{\#}\mid y\times u)x-(x\mid y)x^{\#}\times%
\overline{u} \\
& \quad -(x\mid u)x^{\#}\times\overline{y}+x^{\#}\times\left( x\times \overline{%
(y\times u)}\right) ;
\end{align*}
using (\ref{eq2/11}): $x^{\#}\times\left( x\times\overline{(y\times u)}%
\right) =\det x\ \overline{(y\times u)}+(x^{\#}\mid y\times u)x$, we get
\begin{equation*}
Q(x)D(y,x)u=D(x,y)Q(x)u.
\end{equation*}
This proves (J1).

Let us now prove
\begin{equation}
D(Q(x)y,y)=D(x,Q(y)x).   \tag{J2}
\end{equation}
We have
\begin{align*}
D(Q(x)y,y)z & =(Q(x)y\mid y)z+(z\mid y)Q(x)y-(Q(x)y\times z)\times \overline{%
y} \\
& =((x\mid y)x-x^{\#}\times\overline{y}\mid y)z+(z\mid y)\left( (x\mid
y)x-x^{\#}\times\overline{y}\right) \\
& \quad  -((x\mid y)x-x^{\#}\times\overline{y})\times z)\times\overline{y} \\
& =(x\mid y)^{2}z-2(x^{\#}\mid y^{\#})z+(z\mid y)(x\mid y)x-(z\mid
y)x^{\#}\times\overline{y} \\
& \quad -(x\mid y)\left( x\times z\right) \times\overline{y}+((x^{\#}\times%
\overline{y})\times z)\times\overline{y}
\end{align*}
and
\begin{align*}
D(x,Q(y)x)z & =(x\mid Q(y)x)z+(z\mid Q(y)x)x-(x\times z)\times \overline{%
Q(y)x} \\
& =(x\mid(y\mid x)y-y^{\#}\times\overline{x})z+(z\mid(y\mid x)y-y^{\#}\times%
\overline{x})x \\
& \quad -(x\times z)\times\left( (x\mid y)\overline{y}-\overline{y}^{\#}\times
x\right) \\
& =(x\mid y)^{2}z-2(x^{\#}\mid y^{\#})z+(z\mid y)(x\mid y)x-(z\times x\mid
y^{\#})x \\
& \quad -(x\mid y)\left( x\times z\right) \times\overline{y}+(x\times z)\times(%
\overline{y}^{\#}\times x).
\end{align*}
Applying (\ref{eq2/14}) gives
\begin{equation*}
(x\times z)\times(x\times\overline{y}^{\#})+x^{\#}\times(\overline{y}%
^{\#}\times z)=(x^{\#}\mid y^{\#})z+(x^{\#}\mid z)\overline{y}^{\#}+(z\times
x\mid y^{\#})x;
\end{equation*}
applying (\ref{eq2/14a}) yields
\begin{equation*}
((x^{\#}\times\overline{y})\times z)\times\overline{y}+x^{\#}\times (%
\overline{y}^{\#}\times z)=(x^{\#}\mid y^{\#})z+(x^{\#}\mid z)\overline {y}%
^{\#}+(z\mid y)x^{\#}\times\overline{y}.
\end{equation*}
Comparing these two last identities gives
\begin{equation*}
(x\times z)\times(x\times\overline{y}^{\#})-(z\times x\mid
y^{\#})x=((x^{\#}\times\overline{y})\times z)\times\overline{y}-(z\mid
y)x^{\#}\times\overline{y},
\end{equation*}
which implies $D(Q(x)y,y)z=D(x,Q(y)x)z$ and proves (J2).
\end{proof}

The space $H_{3}(\mathbb{O})$ endowed with the triple product defined by (%
\ref{eq2/18}) will be referred to as the \emph{Hermitian Jordan triple
system }$H_{3}(\mathbb{O)}$, or the \emph{Hermitian JTS of type VI}, or the
\index{exceptional JTS!of dimension 27}\emph{exceptional Hermitian JTS of
dimension }$27.$

The real subspace $H_{3}(\mathbb{O}_{c})$, which is clearly a real Jordan
triple subsystem, and a ``real form'' of $H_{3}(\mathbb{O})$ in the sense
that the triple product in $H_{3}(\mathbb{O})$ is obtained from the product
in $H_{3}(\mathbb{O}_{c})$ by suitable ``complexification'', will be called
the \emph{Euclidean JTS }$H_{3}(\mathbb{O}_{c})$, or the \emph{Euclidean JTS
of type VI}$,$ or the \emph{exceptional compact JTS of dimension }$27$.

\subsection{The minimal polynomial of $H_{3}(\mathbb{O})$}

In this section, we compute the
\index{generic minimal polynomial}\emph{generic minimal polynomial} and the
\index{rank of a Jordan triple}\emph{rank} of the Jordan triple system $H_{3}(\mathbb{O})$ (see
\cite{Roos1999} for the general theory of these notions in a JTS). Recall
that the powers $x^{(k,y)}$ in a Hermitian Jordan triple system $V$ are
defined for $x,y\in V$ and $k\in\mathbb{N}$, $k>0$ by%
\begin{align*}
x^{(1,y)} & =x, \\
 x^{(k+1,y)} & =%
\textstyle\frac{1}{2}D(x,y)x^{(k,y)},
\end{align*}
and the
\index{odd powers in Hermitian JTS}\emph{odd powers} $x^{(2k+1)}$ of $x\in V$, for $k\in%
\mathbb{N}$, by
\begin{equation*}
x^{(2k+1)}=x^{(k+1,x)}.
\end{equation*}
A
\index{tripotent in Hermitian JTS}\emph{tripotent element }in $H_{3}(\mathbb{O})$ is
an element $x$ such that $x^{(3)}=x$.

\begin{lemma}
\label{th2/5}Let $x,y\in H_{3}(\mathbb{O})$. Then
\begin{align}
x^{(2,y)} & =%
\textstyle\frac{1}{2}D(x,y)x=(x\mid y)x-x^{\#}\times\overline {y},  \notag \\
\textstyle\frac{1}{2}D(x,y)\left( x^{\#}\times\overline{y}\right) & =(x^{\#}\mid
y^{\#})x-\det x\ \overline{y}^{\#},  \label{eq2/19} \\
\textstyle\frac{1}{2}D(x,y)\overline{y}^{\#} & =\det\overline{y}\ x.   \label{eq2/20}
\end{align}
The subspace $\sum_{1}^{\infty}\mathbb{C}x^{(k,y)}$ is contained in the
subspace generated by $\left( x,x^{\#}\times\overline{y},\overline{y}%
^{\#}\right) $; the
\index{flat subspace in Hermitian JTS}\emph{flat subspace} generated by $x$:
\begin{equation*}
<<x>>=\sum_{0}^{\infty}\mathbb{C}x^{(2k+1)}
\end{equation*}
is contained in the subspace generated by $\left( x,x^{\#}\times%
\overline {x},\overline{x}^{\#}\right) $.
\end{lemma}

\begin{proof}
The relation for $x^{(2,y)}$ is nothing but the defining relation. From (\ref%
{eq2/18}) and (\ref{eq2/13}), we have
\begin{align*}
D(x,y) & \left( x^{\#}\times\overline{y}\right) =(x\mid y)x^{\#}\times%
\overline{y}+(x^{\#}\times\overline{y}\mid y)x-(x\times(x^{\#}\times%
\overline{y}))\times\overline{y} \\
& =(x\mid y)x^{\#}\times\overline{y}+2(x^{\#}\mid y^{\#})x-((\det x)
\overline{y}+(x\mid y)x^{\#})\times\overline{y} \\
& =2(x^{\#}\mid y^{\#})x-2\det x\ \overline{y}^{\#},
\end{align*}
that is, (\ref{eq2/19}). Using (\ref{eq2/13}) again, we get
\begin{align*}
D(x,y)\overline{y}^{\#} & =(x\mid y)\overline{y}^{\#}+(\overline{y}^{\#}\mid
y)x-(x\times\overline{y}^{\#})\times\overline{y} \\
& =(x\mid y)\overline{y}^{\#}+3\det\overline{y}\ x-\det\overline{y}\
x-(x\mid y)\overline{y}^{\#} \\
& =2\det\overline{y}\ x,
\end{align*}
that is, (\ref{eq2/20}).
\end{proof}

\begin{proposition}
\label{th2/6}The generic minimal polynomial of the Jordan triple system $%
H_{3}(\mathbb{O})$ is
\begin{equation}
m(T,x,y)=T^{3}-(x\mid y)T^{2}+(x^{\#}\mid y^{\#})T-\det x\det\overline{y};
\label{eq2/21}
\end{equation}
the rank of $H_{3}(\mathbb{O})$ is $3$.
\end{proposition}

\begin{proof}
The lemma shows that the JTS $H_{3}(\mathbb{O})$ has rank $\leq3$. It is now
a matter of elementary algebra to compute a linear relation between $x$, $%
x^{(2,y)}$, $x^{(3,y)}$, $x^{(4,y)}$. From (\ref{eq2/19}), (\ref{eq2/20}),
we deduce
\begin{align*}
x^{(3,y)} & =\textstyle\frac{1}{2}D(x,y)x^{(2,y)}=(x\mid y)x^{(2,y)}-(x^{\#}\mid
y^{\#})x+\det x\ \overline{y}^{\#}, \\
x^{(4,y)} & =\textstyle\frac{1}{2}D(x,y)x^{(3,y)}=(x\mid y)x^{(3,y)}-(x^{\#}\mid
y^{\#})x^{(2,y)}+\det x\det\overline{y}\ x.
\end{align*}
This shows that for all $x,y\in V=H_{3}(\mathbb{O})$, the minimal polynomial
of $x$ in $V^{(y)}$ divides
\begin{equation}
T^{3}-(x\mid y)T^{2}+(x^{\#}\mid y^{\#})T-\det x\det\overline{y},
\label{eq2/51}
\end{equation}
and so does the generic minimal polynomial. In order to prove that this is
actually the generic minimal polynomial, we take
\begin{equation*}
x=y=\alpha_{1}e_{1}+\alpha_{2}e_{2}+\alpha_{3}e_{3}
\end{equation*}
with $\alpha_{1}>\alpha_{2}>\alpha_{3}>0$. As it is easily checked, $%
(e_{1},e_{2},e_{3})$ is a set of
\index{orthogonal tripotents}\emph{orthogonal tripotents} (i.e. $%
D(e_{i},e_{i})e_{j}=2\delta_{ij}e_{j}$) and the minimal polynomial of $x$ in
$V^{(x)}$ is $(T-\alpha_{1}^{2})(T-\alpha_{2}^{2})(T-\alpha_{3}^{2})$. This
shows that the generic minimal polynomial has to be of degree $3$ and is
equal to (\ref{eq2/51}).
\end{proof}

\subsection{Positivity; tripotents}

The next proposition implies that $H_{3}(\mathbb{O})$ is a \emph{positive}
Hermitian Jordan triple system.

\begin{proposition}
\label{th2/7}Let $x\in H_{3}(\mathbb{O})$, $x\neq0$. Then $x^{(3)}=\lambda x$
if and only if one of the following occurs:

\begin{enumerate}
\item $(x\mid x)=\lambda$, $x^{\#}=0$;

\item $(x\mid x)=2\lambda$, $(x^{\#}\mid x^{\#})=\lambda^{2}$, $\det x=0$;

\item $(x\mid x)=3\lambda$, $(x^{\#}\mid x^{\#})=3\lambda^{2}$, $\left| \det
x\right| ^{2}=\lambda^{3}.$
\end{enumerate}
\end{proposition}

\begin{proof}
By definition $x^{(3)}=(x\mid x)x-x^{\#}\times
\overline{x}$, so the relation $x^{(3)}=\lambda x$ is equivalent to
\begin{equation}
\left( (x\mid x)-\lambda \right) x=x^{\#}\times \overline{x}.  \label{eq2/22}
\end{equation}%
As $\sum_{1}^{\infty }\mathbb{C}x^{(2k-1)}$ is contained in the subspace
generated by $\left( x,x^{\#}\times \overline{x},\overline{x}^{\#}\right) $,
the relation (\ref{eq2/22}) holds if and only if both sides have the same
Hermitian products with $x$, $x^{\#}\times \overline{x}$, $\overline{x}^{\#}$%
. This provides the conditions
\begin{equation}
\left( (x\mid x)-\lambda \right) (x\mid x)=2(x^{\#}\mid x^{\#}),
\label{eq2/24}
\end{equation}%
\begin{align*}
2\left( (x\mid x)-\lambda \right) (x^{\#}\mid x^{\#})& =(x^{\#}\times
\overline{x}\mid x^{\#}\times \overline{x}) \\
& =3\left\vert \det x\right\vert ^{2}+(x\mid x)(x^{\#}\mid x^{\#}),
\end{align*}%
(using (\ref{eq2/16})), that is,
\begin{equation}
\left( (x\mid x)-2\lambda \right) (x^{\#}\mid x^{\#})=3\left\vert \det
x\right\vert ^{2}  \label{eq2/23}
\end{equation}%
and finally, using
\begin{equation*}
2(x^{\#}\times \overline{x}:x^{\#})=4(\overline{x}:(x^{\#})^{\#})=4(x\mid
x)\det x,
\end{equation*}%
we get%
\begin{equation*}
3\left( (x\mid x)-\lambda \right) \det x=2(x\mid x)\det x,
\end{equation*}%
that is,
\begin{equation}
\left( (x\mid x)-3\lambda \right) \det x=0.  \label{eq2/25}
\end{equation}

If $x\neq0$, $x^{\#}=0$, then (\ref{eq2/23}), (\ref{eq2/25}) are satisfied
and (\ref{eq2/24}) is equivalent to $(x\mid x)=\lambda$.

If $x^{\#}\neq0$ but $\det x=0$, then (\ref{eq2/25}) is satisfied. Condition
(\ref{eq2/23}) is equivalent to $(x\mid x)=2\lambda$ and (\ref{eq2/24}) is
then equivalent to ($x^{\#}\mid x^{\#})=\lambda^{2}$.

If $\det x\neq0$, (\ref{eq2/25}) is equivalent to $(x\mid x)=3\lambda$; (\ref%
{eq2/24}) provides ($x^{\#}\mid x^{\#})=3\lambda^{2}$ and (\ref{eq2/23})
gives $\left| \det x\right| ^{2}=\lambda^{3}$.
\end{proof}

As an immediate consequence, we have

\begin{proposition}
\label{th2/10}The set $\mathcal{E}$ of tripotents of $H_{3}(\mathbb{O})$ is
the disjoint union $\mathcal{E}=\mathcal{E}_{0}\cup\mathcal{E}_{1}\cup%
\mathcal{E}_{2}\cup\mathcal{E}_{3}$, where $\mathcal{E}_{0}=\{0\}$,
\begin{align}
\mathcal{E}_{1} & =\left\{ x\mid(x\mid x)=1,\ x^{\#}=0\right\} ,
\label{eq2/34} \\
\mathcal{E}_{2} & =\left\{ x\mid(x\mid x)=2,\ (x^{\#}\mid x^{\#})=1,\ \det
x=0\right\} ,  \label{eq2/35} \\
\mathcal{E}_{3} & =\left\{ x\mid(x\mid x)=3,\ (x^{\#}\mid x^{\#})=3,\
\left\vert \det x\right\vert ^{2}=1\right\} .   \label{eq2/36}
\end{align}
\end{proposition}

\begin{lemma}
\label{th2/8}Let $x,y\in H_{3}(\mathbb{O})$ be two orthogonal tripotents.
Then $(x\mid y)=0$.
\end{lemma}

\begin{proof}
Let $x,y$ be two tripotents. They are orthogonal if
and only if $D(x,y)=0$.

From (\ref{eq2/20}): $D(x,y)\overline {y}^{\#}=2\det\overline{y}\ x$, we
deduce that $x=0$ or $\det y=0$. So if $y\in\mathcal{E}_{3}$, then $x=0$.
Let $\det y=0$; from (\ref{eq2/19}):
\begin{equation*}
\textstyle\frac{1}{2}D(x,y)\left( x^{\#}\times\overline{y}\right) =(x^{\#}\mid
y^{\#})x-\det x\ \overline{y}^{\#},
\end{equation*}
we then deduce
\begin{equation*}
\textstyle\frac{1}{2}\left( D(x,y)\left( x^{\#}\times\overline{y}\right) \mid y\right)
=(x^{\#}\mid y^{\#})(x\mid y)-3\det x\ \det\overline{y},
\end{equation*}
which implies $(x^{\#}\mid y^{\#})(x\mid y)=0$. Also,
\begin{equation*}
0=x^{(2,y)}=\textstyle\frac{1}{2}D(x,y)x=(x\mid y)x-x^{\#}\times\overline{y},
\end{equation*}
which implies $(x\mid y)x=x^{\#}\times\overline{y}$ and $(x\mid
y)^{2}=2(x^{\#}\mid y^{\#})$; hence
\begin{equation*}
(x\mid y)^{3}=2(x^{\#}\mid y^{\#})(x\mid y)=0.
\end{equation*}
\end{proof}

It follows from Lemma \ref{th2/8} that if $x\in\mathcal{E}_{i}$ and $y\in%
\mathcal{E}_{j}$ are orthogonal tripotents, then their sum $x+y$, which is
also a tripotent, belongs to $\mathcal{E}_{i+j}$. In particular, elements of
$\mathcal{E}_{1}$ are minimal tripotents and elements of $\mathcal{E}_{3}$
are maximal tripotents. The elements of $\mathcal{E}_{3}$ have the following
simple characterization:

\begin{proposition}
\label{th2/9}An element $x\in H_{3}(\mathbb{O})$ is in $\mathcal{E}_{3}$ if
and only if $x\neq0$ and
\begin{equation}
x=\det x\ \overline{x}^{\#}.   \label{eq2/26}
\end{equation}
\end{proposition}

\begin{proof}
If $x\in\mathcal{E}_{3}$, then $2x=x^{\#}\times\overline{x}$, which implies
\begin{equation*}
4x^{\#}=x\times(x^{\#}\times\overline{x})=\det x\ \overline{x}+(x\mid
x)x^{\#}.
\end{equation*}
Hence $x^{\#}=\det x\ \overline{x}$ (as $(x\mid x)=3$) and $\det\overline {x}%
\ x^{\#}=\overline{x}$ (as $\left| \det x\right| ^{2}=1).$

Conversely, let $x\neq0$, $x=\det x\ \overline{x}^{\#}$; this implies $\det
x\neq0$. Then $x^{\#}\times\overline{x}=x^{\#}\times\left( \det\overline {x}%
\ x^{\#}\right) =2\left| \det x\right| ^{2}x$. This means that $%
x^{(3)}=\lambda x$, with $\lambda=(x\mid x)-2\left| \det x\right| ^{2}$. By
Proposition \ref{th2/7}, we have $\left| \det x\right| ^{2}=\lambda^{3}$ and
$(x\mid x)=3\lambda$, which implies $(x\mid x)=3\left| \det x\right| ^{2}$
and $\lambda=\lambda^{3}>0$; hence $\lambda=1$ and $x\in\mathcal{E}_{3}$.
\end{proof}

For a tripotent $x$ in a positive Hermitian JTS $V$, the operator $D(x,x)$ is self-adjoint
and its eigenvalues are in $\{0, 1, 2\}$. The
\index{PeirceC@\textsc{Peirce, Charles Sanders} (1839--1914)}%
\index{Peirce decomposition}%
\emph{Peirce decomposition} relative to $x$ is
\[
V=V_{2}(x)\oplus V_{1}(x)\oplus V_{0}(x),
\]
where the \emph{Peirce subspaces} $V_{j}(x)$ are the eigenspaces of $D(x,x)$:%
\[
V_{j}(x)=\left\{  y\in V\mid D(x,x)y=jy\right\}  ,\quad j=0,1,2.
\]

\begin{proposition}
\label{th2/11}For $x\in\mathcal{E}_{3}$, the Peirce subspaces are $%
V_{0}(x)=V_{1}(x)=0$, $V_{2}(x)=V=H_{3}(\mathbb{O}).$
\end{proposition}

\begin{proof}
It suffices to prove that for each $y\in H_{3}(\mathbb{O}),$ one has $%
D(x,x)y=2y$. As $D(x,x)x=2x$, it is enough to prove this if $(x\mid y)=0$.
If $(x\mid y)=0$, we have then, using (\ref{eq2/26}) and (\ref{eq2/11}),
\begin{align*}
D(x,x)y & =(x\mid x)y-(x\times y)\times\overline{x} \\
& =3y-\det\overline{x}\ (x\times y)\times x^{\#} \\
& =3y-\det\overline{x}\ \left( \det x\ y+(x^{\#}:y)x\right) =2y,
\end{align*}
as $\left| \det x\right| ^{2}=1$ and, by (\ref{eq2/26}), $\det\overline {x}\
(x^{\#}:y)=(y\mid x)=0$.
\end{proof}

Proposition \ref{th2/6} shows that a maximal flat subspace has dimension $3$%
. From Proposition \ref{th2/10}, we see that $e_{1}$, $e_{2}$, $e_{3}$
belong to $\mathcal{E}_{1}$ and are therefore minimal tripotents. From the
definition
\begin{equation*}
D(e_{1},e_{1})z=z+(z\mid e_{1})e_{1}-e_{1}\times(e_{1}\times z)
\end{equation*}
and from the relations (\ref{eq2/27}), it is easily checked that
\begin{align}
V_{0}(e_{1}) & =\mathbb{C}e_{2}\oplus\mathbb{C}e_{3}\oplus\mathcal{F}_{1},
\label{eq2/37} \\
V_{1}(e_{1}) & =\mathcal{F}_{2}\oplus\mathcal{F}_{3},  \label{eq2/38} \\
V_{2}(e_{1}) & =\mathbb{C}e_{1}.   \label{eq2/39}
\end{align}
Similar results hold for the Peirce decomposition with respect to $e_{2}$
and $e_{3}$. As $e_{2}$ and $e_{3}$ belong to $V_{0}(e_{1})$, they are
orthogonal to $e_{1}$; also, $e_{2}$ is orthogonal to $e_{3}$. So $%
(e_{1},e_{2},e_{3})$ is a frame for the Jordan triple system $H_{3}(\mathbb{O%
})$ and $\mathbb{R}e_{1}\oplus\mathbb{R}e_{2}\oplus\mathbb{R}e_{3}$ is a
maximal flat subspace. It is also easily checked that the simultaneous
Peirce decomposition with respect to the frame $(e_{1},e_{2},e_{3})$ is
\begin{equation*}
H_{3}(\mathbb{O)=}\bigoplus_{1\leq i\leq j\leq3}V_{ij},
\end{equation*}
with $V_{ii}=\mathbb{C}e_{i}$, $V_{ij}=\mathcal{F}_{k}$.

\begin{theorem}[Freudenthal's theorem]
\label{th2/12}\index{Freudenthal's theorem}%
\index{FreudenthalH@\textsc{Freudenthal, Hans} (1905--1990)}%
Let $x\in H_{3}(\mathbb{O})$. Then there exists $k\in\operatorname{Aut}%
H_{3}(\mathbb{O})$ such that
\begin{equation*}
kx=\alpha_{1}e_{1}+\alpha_{2}e_{2}+\alpha_{3}e_{3}\quad(\alpha_{1},\alpha
_{2},\alpha_{3}\in\mathbb{R}).
\end{equation*}
\end{theorem}

Actually, $H_{3}(\mathbb{O})$ is a positive Jordan triple system; then for
each $x$ there exists an automorphism $k$ such that $kx$ belongs to the
maximal flat subspace $\mathbb{R}e_{1}\oplus\mathbb{R}e_{2}\oplus \mathbb{R}%
e_{3}$. The theorem may also be proved directly in this special case,
following the lines of the general theory.

\begin{theorem}
\label{th2/13}The Hermitian Jordan triple system $H_{3}(\mathbb{O})$ is
simple and of tube type. Its numerical invariants are
\begin{equation*}
a=8,\ b=0,\ r=3,\ g=18.
\end{equation*}
\end{theorem}

To show that $H_{3}(\mathbb{O})$ is simple, it is enough to find a frame
such that all $V_{ij}$ $(1\leq i<j\leq3)$ are non-zero; this occurs with the
frame $(e_{1},e_{2},e_{3})$. We then have $a=\dim\mathcal{F}_{i}=8$, $b=\dim
V_{0i}=0$, $g=2+a(r-1)=18$.

\begin{corollary}
In the Jordan triple system $H_{3}(\mathbb{O})$, we have%
\begin{align}
\operatorname{Tr}D(x,y) & =18(x\mid y),  \label{eq2/33} \\
\operatorname{Det}B(x,y) & =\left( 1-(x\mid y)+(x^{\#}\mid y^{\#})-\det x\det%
\overline{y}\right) ^{18},   \label{eq2/32}
\end{align}
where $\operatorname{Tr}$ and $\operatorname{Det}$ denote the trace and determinant of $%
\mathbb{C}$-linear operators in $H_{3}(\mathbb{O})$.
\end{corollary}

\subsection{The exceptional Jordan triple system of dimension $16$\label%
{Exc16}}

We consider the subsystem of $H_{3}(\mathbb{O)}$%
\begin{equation*}
V_{1}(e_{1})=\mathcal{F}_{2}\oplus\mathcal{F}_{3},
\end{equation*}
which is then a positive Hermitian Jordan triple system of dimension $16$.
Let us denote the space $V_{1}(e_{1})$ by $W$. For $%
x=F_{2}(x_{2})+F_{3}(x_{3})\in W$, we have, according to (\ref{eq2/4}),
\begin{equation*}
x^{\#}=-n(x_{2})e_{2}-n(x_{3})e_{3}+\widetilde{F_{1}}(x_{2}x_{3})\in
V_{0}(e_{1})
\end{equation*}
and $\det x=0$. The structure of Jordan triple system in $W$ is defined, for
$x=F_{2}(x_{2})+F_{3}(x_{3})$, $y=F_{2}(y_{2})+F_{3}(y_{3})$, by
\begin{align*}
Q(x)y & =(x\mid y)x-x^{\#}\times\overline{y} \\
& =\left( (x_{2}\mid y_{2})+(x_{3}\mid y_{3})\right) x-n(x_{2})F_{2}(%
\overline{y_{2}})-n(x_{3})F_{3}(\overline{y_{3}}) \\
& \quad -\widetilde{F_{3}}((\widetilde{x_{3}}\widetilde{x_{2}})\overline{y_{2}})-%
\widetilde{F_{2}}(\overline{y_{3}}(\widetilde{x_{3}}\widetilde{x_{2}})) \\
& =F_{2}\left( (x_{2}\mid y_{2})x_{2}+(x_{3}\mid y_{3})x_{2}-n(x_{2})%
\overline{y_{2}}-(x_{2}x_{3})\widetilde{\overline{y_{3}}}\right) \\
& \quad +F_{3}\left( (x_{2}\mid y_{2})x_{3}+(x_{3}\mid y_{3})x_{3}-n(x_{3})%
\overline{y_{3}}-\widetilde{\overline{y_{2}}}(x_{2}x_{3})\right) .
\end{align*}
Using identities in Cayley algebras, we get%
\begin{equation}
Q(x)y=F_{2}\left( x_{2}\widetilde{\overline{y_{2}}}x_{2}+(x_{2}\overline{%
y_{3}})\widetilde{x_{3}}\right) +F_{3}\left( \widetilde{x_{2}}(\overline{%
y_{2}}x_{3})+x_{3}\widetilde{\overline{y_{3}}}x_{3}\right) .   \label{eq2/28}
\end{equation}
The triple product in $W$ is then given by
\begin{align}
\{xyz\} & =F_{2}\left( (x_{2}\widetilde{\overline{y_{2}}})z_{2}+(z_{2}%
\widetilde{\overline{y_{2}}})x_{2}+(x_{2}\overline{y_{3}})\widetilde{z_{3}}%
+(z_{2}\overline{y_{3}})\widetilde{x_{3}}\right) \notag \\
& \quad +F_{3}\left( \widetilde{x_{2}}(\overline{y_{2}}z_{3})+\widetilde{z_{2}}(%
\overline{y_{2}}x_{3})+x_{3}(\widetilde{\overline{y_{3}}}z_{3})+z_{3}(%
\widetilde{\overline{y_{3}}}x_{3})\right) .   \label{eq2/30}
\end{align}

\begin{proposition}
\label{th2/14}The generic minimal polynomial of $W$ is
\begin{equation*}
m_{W}(T;x,y)=T^{2}-(x\mid y)T+(x^{\#}\mid y^{\#}).
\end{equation*}
For $x,y\in W$, the subspace $\sum_{1}^{\infty}\mathbb{C}x^{(k,y)}$ is
contained in $\mathbb{C}x+\mathbb{C}x^{\#}\times\overline{y}$; the flat
subspace generated by $x$:
\begin{equation*}
<<x>>=\sum_{0}^{\infty}\mathbb{C}x^{(2k+1)}
\end{equation*}
is contained in the subspace generated by $\left( x,x^{\#}\times\overline {x}%
\right) $.
\end{proposition}

\begin{proof}
Let $x,y\in W$; then, by (\ref{eq2/19}) and $\det x=0$, we have
\begin{align*}
& x^{(2,y)}  =\textstyle\frac{1}{2}D(x,y)x=(x\mid y)x-x^{\#}\times\overline{y}, \\
& \textstyle\frac{1}{2}D(x,y)\left( x^{\#}\times\overline{y}\right) =(x^{\#}\mid
y^{\#})x.
\end{align*}
This shows $\sum_{1}^{\infty}\mathbb{C}x^{(k,y)}\subset\mathbb{C}x+\mathbb{C}%
x^{\#}\times\overline{y}.$ Moreover, these relations imply
\begin{equation*}
x^{(3,y)}=(x\mid y)x^{(2,y)}-(x^{\#}\mid y^{\#})x,
\end{equation*}
which shows that the generic minimal polynomial $m(T,x,y)=m_{W}(T;x,y)$
divides $T^{2}-(x\mid y)T+(x^{\#}\mid y^{\#})$. To prove equality, it will
be enough to prove that the rank of $W$ is $2$, that is, to find $x\in W$
such that $x$ and $x^{\#}\times\overline{x}$ are $\mathbb{R}$-linearly
independent. For this, take $x=F_{2}(b)$, with $b\in\mathbb{O}$; then $%
x^{\#}=-n(b)e_{2}$ and $x^{\#}\times\overline{x}=n(b)F_{2}(\overline{b})$.
So it suffices to choose $b\in\mathbb{O}$ such that $n(b)=1$ and $b,%
\overline{b}$ linearly independent. This proves $m(T,x,y)=T^{2}-(x\mid
y)T+(x^{\#}\mid y^{\#})$.
\end{proof}

A tripotent in $W$ is also a tripotent in $V=H_{3}(\mathbb{O)}$; as $\det x=0
$ for each $x\in W$, $\mathcal{E}_{3}\cap W=\emptyset$, the set of
tripotents of $W$ is $\mathcal{E}^{\prime}=\mathcal{E}_{0}^{\prime}\cup%
\mathcal{E}_{1}^{\prime}\cup\mathcal{E}_{2}^{\prime}$ with $\mathcal{E}%
_{j}^{\prime }=\mathcal{E}_{j}\cap W$. Also, two orthogonal tripotents $%
x,y\in W$ are orthogonal in $V$ and hence verify $(x\mid y)=0$; it follows
that elements of $\mathcal{E}_{1}^{\prime}$ are minimal tripotents and
elements of $\mathcal{E}_{2}^{\prime}$ are maximal tripotents.

Minimal tripotents $F_{2}(\beta)+F_{3}(\gamma)$ are characterized by
\begin{equation}
n(\beta)=n(\gamma)=0,\quad\beta\gamma=0,\quad(\beta\mid\beta)+(\gamma
\mid\gamma)=1.   \label{eq2/40}
\end{equation}

An example of minimal tripotent is then given by $u=F_{2}(\beta)$, with $%
\beta$ satisfying
\begin{equation}
(\beta\mid\beta)=1,\quad n(\beta)=0;   \label{eq2/41}
\end{equation}
these relations are equivalent to
\begin{equation}
\textstyle \beta=b_{1}+\mathrm{i}b_{2},\quad b_{1},b_{2}\in\mathbb{%
O}_{c},\quad n(b_{1})=n(b_{2})=\frac{1}{4},\quad(b_{1}:b_{2})=0.
\label{eq2/42}
\end{equation}

\begin{lemma}
\label{th2/15}Let $\beta\in\mathbb{O}$ such that $(\beta\mid\beta)=1$ and $%
n(\beta)=0$. Then for each $x\in\mathbb{O}$,
\begin{equation}
x=\widetilde{\beta}(\overline{\beta}x)+\widetilde{\overline{\beta}}(\beta
x).   \label{eq2/29}
\end{equation}
If $L(\beta)$ denotes the left multiplication by $\beta$ in $\mathbb{O}$: $%
L(\beta)x=\beta x$, the following direct sum decomposition holds:
\begin{equation}
\mathbb{O}=\ker L(\beta)\oplus\ker L(\overline{\beta});   \label{eq2/31}
\end{equation}
moreover,
\begin{equation*}
\ker L(\beta)=\operatorname{Im}L(\widetilde{\beta})
\end{equation*}
and
\begin{equation*}
\dim_{\mathbb{C}}\ker L(\beta)=4.
\end{equation*}
\end{lemma}

\begin{proof}
Polarizing the identity $n(y)x=\widetilde{y}(yx)$ in $\mathbb{O}$, we obtain
\begin{equation*}
(y:z)x=\widetilde{y}(zx)+\widetilde{z}(yx)
\end{equation*}
and (\ref{eq2/29}) follows by $y\leftarrow\beta$, $z\leftarrow\overline{%
\beta }.$ Let $x_{1}=\widetilde{\beta}(\overline{\beta}x)$ and $x_{2}=%
\widetilde {\overline{\beta}}(\beta x)$; then $x_{1}\in\ker L(\beta)$ and $%
x_{2}\in\ker L(\overline{\beta})$. If $x\in\ker L(\beta)\cap\ker L(\overline{%
\beta})$, it follows from (\ref{eq2/29}) that $x=0$; this proves (\ref%
{eq2/31}).

Clearly $\operatorname{Im}L(\widetilde{\beta})\subset\ker L\left( \beta\right) $, by
$\beta(\widetilde{\beta}x)=n\left( \beta\right) x=0$. Assume $\beta x=0$; by
(\ref{eq2/29}), we have $x=\widetilde{\beta}(\overline{\beta}x)$, that is $%
x\in\operatorname{Im}L(\widetilde{\beta})$. So $\ker L(\beta)=\operatorname{Im}L(\widetilde{%
\beta})$.

As $x\mapsto\overline{x}$ is a (real) automorphism of $\mathbb{O}$, the
spaces $\ker L(\beta)$ and $\ker L(\overline{\beta})$ have the same real
dimension, hence also the same complex dimension. This implies $\dim_{%
\mathbb{C}}\ker L(\beta)=4$.
\end{proof}

We are now able to compute the Peirce decomposition with respect to the
minimal tripotent $u=F_{2}(\beta)$.

\begin{lemma}
\label{th2/16}Let $\beta\in\mathbb{O}$ such that $(\beta\mid\beta)=1$ and $%
n(\beta)=0$. The spaces of the Peirce decomposition of $W$ with respect to $%
u=F_{2}(\beta)$ are
\begin{align}
W_{0}(u) & =\mathbb{C}\overline{u}\oplus F_{3}(\ker L( \overline {\beta}%
) ),  \label{eq2/43} \\
W_{1}(u) & =F_{2}(\langle \beta,\overline{\beta}\rangle ^{\perp
})\oplus F_{3}(\ker L\left( \beta\right) ),  \label{eq2/44} \\
W_{2}(u) & =\mathbb{C}u.   \label{eq2/45}
\end{align}
\end{lemma}

Here $\left\langle \beta,\overline{\beta}\right\rangle ^{\perp}$ stands for
the orthogonal subspace of $\mathbb{O}$, with respect to the Hermitian
product $\left( ~\mid~\right) $, of the $2$-dimensional subspace $\mathbb{%
C\beta \oplus C}\overline{\beta}$. Note that the conditions on $\beta$ mean
that $\left( \beta,\overline{\beta}\right) $ is orthonormal.

\begin{proof}
For $x=F_{2}(x_{2})+F_{3}(x_{3})$, we have
\begin{align*}
D(u,u)x & =(u\mid u)x+(x\mid u)u-(u\times x)\times\overline{u} \\
 =\ & F_{2}(x_{2})+F_{3}(x_{3})+(x_{2}\mid\beta)F_{2}(\beta) \\
& \quad -\left( F_{2}(\beta)\times F_{2}(x_{2})\right) \times F_{2}(\overline{\beta%
})-\left( F_{2}(\beta)\times F_{3}(x_{3})\right) \times F_{2}(\overline{\beta%
}) \\
 =\ & F_{2}(x_{2})  +F_{3}(x_{3})+(x_{2}\mid\beta)F_{2}(\beta)+(x_{2}:\beta
)F_{2}(\overline{\beta})-F_{3}(\widetilde{\overline{\beta}}(\beta x_{3}));
\end{align*}
finally
\begin{align*}
D(u,u)x & =F_{2}\left( x_{2}+(x_{2}\mid\beta)\beta-(x_{2}:\beta )\overline{%
\beta}\right) +F_{3}(x_{3}-\widetilde{\overline{\beta}}(\beta x_{3})) \\
& =F_{2}\left( x_{2}+(x_{2}\mid\beta)\beta-(x_{2}:\beta)\overline{\beta }%
\right) +F_{3}(\widetilde{\beta}(\overline{\beta}x_{3})).
\end{align*}
From these two expressions of $D(u,u)x$, it is easily seen that

\begin{itemize}
\item $D(u,u)x=0$ if $x_{2}\in\mathbb{C}\overline{\beta}$ and $\overline {%
\beta}x_{3}=0$;

\item $D(u,u)x=x$ if $(x_{2}\mid\beta)=(x_{2}\mid\overline{\beta})=0$ and $%
\beta x_{3}=0$;

\item $D(u,u)x=2x$ if $x_{2}\in\mathbb{C}\beta$ and $x_{3}=0$.
\end{itemize}

This provides the diagonalization of $D(u,u)$ with the indicated eigenspaces.
\end{proof}

It follows easily from Lemma \ref{th2/16} that $v=F_{2}(\overline{\beta})$
is a minimal tripotent, orthogonal to $u$. The eigenspaces of $D(v,v)$ are
obtained from Lemma \ref{th2/16} with $\beta\leftarrow\overline{\beta}$. By
comparing the two Peirce decompositions, we obtain

\begin{proposition}
\label{th2/17}The spaces of the simultaneous Peirce decomposition with
respect to the frame $\left( u,v\right) =\left( F_{2}(\beta),F_{2}(\overline{%
\beta })\right) $ are
\begin{align}
W_{01} & =F_{3}(\ker L(\beta)),\quad W_{02}=F_{3}(\ker L\left( \overline{%
\beta}\right) ),  \label{eq2/46} \\
W_{12} & =F_{2}(\left\langle \beta,\overline{\beta}\right\rangle ^{\perp
}),\quad W_{11}=\mathbb{C}F_{2}(\beta),\quad W_{22}=\mathbb{C}F_{2}(%
\overline{\beta}).   \label{eq2/47}
\end{align}
\end{proposition}

\begin{proposition}
\label{th2/18}The triple system $W$ is simple. Its numerical invariants are $%
a=6$, $b=4$, $r=2$, $g=12$. In $W$,
\begin{align*}
\operatorname{Tr}D(x,y) & =12(x\mid y), \\
\operatorname{Det}B(x,y) & =\left( 1-(x\mid y)+(x^{\#}\mid y^{\#})\right) ^{12}.
\end{align*}
The set of tripotents of $W$ is $\mathcal{E}^{\prime}=\mathcal{E}%
_{0}^{\prime }\cup\mathcal{E}_{1}^{\prime}\cup\mathcal{E}_{2}^{\prime}$,
with $\mathcal{E}_{0}^{\prime}=\{0\}$,
\begin{align}
\mathcal{E}_{1}^{\prime} & =\left\{ x\in W;\ (x\mid x)=1,\ x^{\#}=0\right\} ,
\label{eq2/48} \\
\mathcal{E}_{2}^{\prime} & =\left\{ x\in W;\ (x\mid x)=2,\ (x^{\#}\mid
x^{\#})=1\right\} .   \label{eq2/49}
\end{align}
\end{proposition}

\begin{proof}
The tripotents of $W$ have already been described. From the previous
proposition, we see that $\dim W_{12}=6$ for the frame $\left(
F_{2}(\beta),F_{2}(\overline{\beta})\right) $. This implies that $W$ is
simple, as it is positive as a subsystem of the positive Hermitian JTS $%
H_{3}(\mathbb{O)}$. The numerical invariants are $r=2$, $a=\dim W_{12}=6$, $%
b=\dim W_{01}=4$, $g=2+a(r-1)+b=12$.
\end{proof}

As an example of maximal tripotent of $W$, we have $w=u+v=F_{2}(\beta +%
\overline{\beta})=F_{2}(c)$, with $c\in\mathbb{O}_{c}$ and $n(c)=1$. The
Peirce spaces for $w$ are $W_{2}(w)=\mathcal{F}_{2}$, $W_{1}(w)=\mathcal{F}%
_{3}$.

The simple positive JTS $W$ is called the
\index{exceptional JTS!of dimension 16}$\emph{exceptional}$\emph{\ JTS of
dimension }$16$.\medskip

Let us look at the Jordan structure of the Peirce subspaces with respect to
the minimal tripotent $u$. The subspace $W_{0}(u)$ has rank $1$ and is
isomorphic to $I_{1,5}$. Consider
\begin{equation*}
W^{\prime}=W_{1}(u)=F_{2}(\left\langle \beta,%
\overline{\beta}\right\rangle ^{\perp})\oplus F_{3}(\ker L\left(
\beta\right) ).
\end{equation*}

Let $\gamma\in\left\langle \beta,\overline{\beta}\right\rangle ^{\perp}$such
that $n(\gamma)=0$ and $\left( \gamma\mid\gamma\right) =1$. Then $u^{\prime
}=F_{2}(\gamma)$ and $v^{\prime}=F_{2}(\overline{\gamma})$ are two
orthogonal tripotents in $W^{\prime}$ and form a frame for $W^{\prime}$. The
spaces of the total Peirce decomposition of $W^{\prime}$ with respect to
this frame are obtained from the corresponding spaces in $W$ by intersection
with $W^{\prime}$, which gives%
\begin{align*}
W_{01}^{\prime} & =F_{3}(\ker L(\beta)\cap\ker L(\gamma)),\quad
W_{02}^{\prime}=F_{3}(\ker L\left( \beta\right) \cap\ker L(\overline{\gamma }%
)), \\
W_{12}^{\prime} & =F_{2}(\left\langle \beta,\overline{\beta},\gamma ,%
\overline{\gamma}\right\rangle ^{\perp}),\quad W_{11}^{\prime}=\mathbb{C}%
F_{2}(\gamma),\quad W_{22}=\mathbb{C}F_{2}(\overline{\gamma}).
\end{align*}
Clearly $\dim W_{12}^{\prime}=4$, which implies that $W^{\prime}$ is simple;
then $\dim W_{01}^{\prime}=\dim W_{02}^{\prime}=2$.

The only simple positive Hermitian Jordan triple system with rank $2$,
dimension $10$ and $a=4$ is $II_{5}$. This proves

\begin{proposition}
The Peirce subspace $W_{1}(u)$ of the exceptional JTS of type $V$ with
respect to a minimal tripotent $u$ is of type $II_{5}$.
\end{proposition}

\begin{exercise}
Prove this proposition directly.
\end{exercise}

\begin{exercise}
Consider $v^{\prime}=F_{3}(\gamma)$ with $\gamma$ subject to the same
conditions as $\beta$:
\begin{equation*}
(\gamma\mid\gamma)=1,\quad n(\gamma)=0.
\end{equation*}
Then $v^{\prime}$ is another minimal tripotent.

\begin{enumerate}
\item Show that $v^{\prime}$ is orthogonal to $u=F_{2}(\beta)$ if and only
if $\overline{\beta}\gamma=0$.

\item Compute the simultaneous Peirce decomposition with respect to the
frame $(F_{2}(\beta),\widetilde{F_{3}}(\overline{\beta})).$
\end{enumerate}
\end{exercise}

\begin{exercise}
Let $\beta\in\mathbb{O}$ such that $(\beta\mid\beta)=1$ and $n(\beta)=0$.

\begin{enumerate}
\item Compute the Peirce decomposition of the JTS of type VI $H_{3}(\mathbb{O%
})$ with respect to the minimal tripotent $u=F_{2}(\beta)$.

\item Find a minimal tripotent $f$ such that $\left( F_{2}(\beta ),F_{2}(%
\overline{\beta}),f\right) $ is a frame of $H_{3}(\mathbb{O})$.
\end{enumerate}
\end{exercise}

\begin{exercise}
Find a Jordan triple subsystem $W^{\prime}$ of $H_{3}(\mathbb{O})$,
isomorphic to $W$, containing $e_{1}$ and $e_{2}$. Compute the Peirce
decomposition of $W^{\prime}$ with respect to $\left( e_{1},e_{2}\right) $.
\end{exercise}


\section{The exceptional symmetric domains}



\subsection{Description of exceptional symmetric domains}

We apply the general results of \cite{Loos1977}. As in the previous section,
we denote by $\mathbb{O}$ the algebra of complex octonions, by $V=H_{3}(%
\mathbb{O})$ the exceptional Jordan system with the Jordan triple structure
defined by Definition \ref{def2/5}, by $W=\mathcal{F}_{2}\oplus\mathcal{F}%
_{3}$ the subsystem of dimension $16$ studied in Subsection \ref{Exc16}.
Recall that these two complex Jordan triples are Hermitian positive and
simple, with respective generic minimal polynomials
\begin{align}
m_{V}(T,x,y) & =T^{3}-(x\mid y)T^{2}+(x^{\#}\mid y^{\#})T-\det x\det\overline{y},
\label{eq3/1}\\
m_{W}(T,x,y) & =T^{2}-(x\mid y)T+(x^{\#}\mid y^{\#}).   \label{eq3/2}
\end{align}

For a Hermitian positive Jordan triple of rank $r$ and generic minimal
polynomial $m(T,x,y)$, the associated circled bounded symmetric domain is
defined by the $r$ inequalities
\begin{equation}
f_{k+1}(x,x)\equiv\frac{1}{k!}\left. \frac{\operatorname{d}^{k}}{\operatorname{d}T^{k}}%
m(T;x,x)\right\vert _{T=1}>0\qquad(k=0,\ldots,r-1).   \label{eq3/3}
\end{equation}

It follows that the symmetric domain $\Omega=\Omega_{V}$ associated to $V$
(called the
\index{exceptional symmetric domain !of dimension 27}\emph{exceptional
symmetric domain of dimension }$27$, or the \emph{symmetric domain of type }$%
VI$) is the set of points in $H_{3}(\mathbb{O})$ which satisfy
\begin{align}
f_{1}(x,x) & \equiv1-(x\mid x)+(x^{\#}\mid x^{\#})-\left\vert \det
x\right\vert ^{2}>0,  \label{ineq27-1} \\
f_{2}(x,x) & \equiv3-2(x\mid x)+(x^{\#}\mid x^{\#})>0,  \label{ineq27-2} \\
f_{3}(x,x) & \equiv3-(x\mid x)>0,   \label{ineq27-3}
\end{align}
while the symmetric domain $\Omega^{\prime}=\Omega_{W}$ associated to $W$
(called the
\index{exceptional symmetric domain !of dimension 16}\emph{exceptional
symmetric domain of dimension }$16$, or the \emph{symmetric domain of type }$%
V$) is the set of points in $W=\mathcal{F}_{2}\oplus \mathcal{F}_{3}$ which
satisfy
\begin{align}
g_{1}(x,x) & \equiv1-(x\mid x)+(x^{\#}\mid x^{\#})>0,  \label{ineq16-1} \\
g_{2}(x,x) & \equiv2-(x\mid x)>0.   \label{ineq16-2}
\end{align}

\subsection{Structure of the boundary}

\subsubsection{General results}

The inequalities (\ref{eq3/3}) are equivalent to the fact that all roots of
the polynomial $m(T;x,x)$ in $T$ (which are always positive) are less than $1
$. The boundary of the symmetric domain $\Omega$ is the disjoint union of
locally closed submanifolds $\partial_{k}\Omega$, which correspond to the
case where $1$ is a root of $m(T;x,x)$ with multiplicity $k$ and the
remaining roots are less than $1$. We first recall general results, valid
for each simple Hermitian positive JTS and the associated irreducible
bounded symmetric domain (see for example \cite{Loos1977}, \S \S 5-6). Then
we apply these results to the case of the two exceptional symmetric domains.

The description of the boundary (Proposition \ref{prop3/8}) also involves
the manifold of tripotents of the corresponding JTS, which is described in
Proposition \ref{prop3/7} below. The description of their tangent space
needs a refinement of the
\index{Peirce decomposition}Peirce decomposition
$V=V_{2}(e)\oplus V_{1}(e)\oplus V_{0}(e)$
associated to a tripotent $e$.

\begin{proposition}
\label{prop3/6}For a tripotent $e$ in a Hermitian positive Jordan triple
system $V$, the operator $Q(e)$ is zero on $V_{1}(e)\oplus V_{0}(e)$ and
restricts to a $\mathbb{C}$-antilinear involution on $V_{2}(e)$.
\end{proposition}

Denoting by $V_{2}^{+}(e)$, $V_{2}^{-}(e)$ the eigenspaces of $Q(e)$ for the
eigenvalues $+1$, $-1$:%
\begin{align*}
V_{2}^{+}(e) & =\left\{ x\in V\mid Q(e)x=x\right\} , \\
V_{2}^{-}(e) & =\left\{ x\in V\mid Q(e)x=-x\right\} ,
\end{align*}
we have $V_{2}^{-}(e)=%
\operatorname{i}V_{2}^{+}(e)$ and $V_{2}(e)=V_{2}^{+}(e)\oplus V_{2}^{-}(e)$. In any
simple Hermitian positive Jordan triple system, one then has (see \cite%
{Loos1977}, Theorem 5.6):

\begin{proposition}
\label{prop3/7}The set $\mathcal{E}_{k}$ of tripotents of rank $k$ is a
compact connected submanifold, and the group $K$ of automorphisms of the
Jordan triple system acts transitively on $\mathcal{E}_{k}$. For $e\in%
\mathcal{E}_{k}$, the direction of the tangent space to $\mathcal{E}_{k}$ at
$e$ is
\begin{equation*}
\overrightarrow{T_{e}\mathcal{E}_{k}}=\operatorname{i}V_{2}^{+}(e)\oplus V_{1}(e).
\end{equation*}
The complex tangent space $H_{e}\mathcal{E}_{k}$ to $\mathcal{E}_{k}$ at $e$
has direction
\begin{equation*}
\overrightarrow{H_{e}\mathcal{E}_{k}}=V_{1}(e).
\end{equation*}
The manifold $\mathcal{E}_{k}$ is a Cauchy-Riemann manifold of CR type $(s,t)
$ and real dimension $2s+t$, with
\begin{equation*}
s=\dim_{\mathbb{C}}V_{1}(e),\quad t=\dim_{\mathbb{R}}V_{2}^{+}(e)=\dim _{%
\mathbb{C}}V_{2}(e).
\end{equation*}
\end{proposition}

\begin{proposition}
\label{prop3/8}Let $V$ be a simple Hermitian positive JTS of rank $r$ and
let $\Omega$ be the associated irreducible symmetric domain. The boundary $%
\partial\Omega$ of the symmetric domain $\Omega=\Omega_{V}$ of type $VI$ is
the disjoint union%
\begin{equation}
\partial\Omega={\displaystyle\bigsqcup\limits_{k=1}^{r}} \partial_{k}\Omega
\label{strat}
\end{equation}
of locally closed manifolds%
\begin{equation*}
\partial_{k}\Omega=\left\{ x\in V\mid f_{j}(x,x)=0\quad(1\leq j\leq k),\
f_{m}(x,x)>0\quad(m>k)\right\} ,
\end{equation*}
where the $f_{j}$'s are defined by (\ref{eq3/3}). The \textquotedblleft
boundary part\textquotedblright\ $\partial_{k}\Omega$ contains the manifold $%
\mathcal{E}_{k}$ of rank $k$ tripotents in $V$. Each $\partial_{k+1}\Omega$
is contained in $\overline{\partial_{k}\Omega}\setminus\partial_{k}\Omega$.

For $e\in\mathcal{E}_{k}$, the normal direction at $e$ to $\partial_{k}\Omega
$ is $V_{2}^{+}(e)$ and the direction $\overrightarrow{T_{e}(\partial_{k}%
\Omega)}$ of the tangent space $T_{e}(\partial_{k}\Omega)$ is
\begin{equation*}
\overrightarrow{T_{e}(\partial_{k}\Omega)}=\operatorname{i}V_{2}^{+}(e)\oplus
V_{1}(e)\oplus V_{0}(e).
\end{equation*}

The intersection of $\partial_{k}\Omega$ with the affine tangent space $e+%
\overrightarrow{T_{e}(\partial_{k}\Omega)}$ is
\begin{equation*}
\partial_{k}\Omega\cap\left( e+\overrightarrow{T_{e}(\partial_{k}\Omega )}%
\right) =\partial_{k}\Omega\cap\left( e+V_{0}(e)\right) =e+\Omega(e),
\end{equation*}
where $\Omega(e)$ is the symmetric domain associated to the Jordan triple
subsystem $V_{0}(e)$. The direction of the tangent space to $\partial
_{k}\Omega$ is constant along $e+\Omega(e)$.

The boundary part $\partial_{k}\Omega$ is the disjoint union
\begin{equation*}
\partial_{k}\Omega={\displaystyle\bigsqcup\limits_{e\in\mathcal{E}_{k}}}
\left( e+\Omega(e)\right) .
\end{equation*}
Let
\begin{equation*}
p_{k}:\partial_{k}\Omega\rightarrow\mathcal{E}_{k}
\end{equation*}
be defined by $p_{k}(x)=e$ if $e\in\mathcal{E}_{k}$ and $x-e\in V_{0}(e)$.
Then $\left( \partial_{k}\Omega,\mathcal{E}_{k},p_{k}\right) $ is a locally
trivial fiber bundle, isomorphic to $\left( \mathcal{X}_{k},\mathcal{E}%
_{k},q_{k}\right) $, where
\begin{equation*}
\mathcal{X}_{k}=\left\{ (e,y)\in\mathcal{E}_{k}\times V\mid y\in
\Omega(e)\right\}
\end{equation*}
and $q_{k}$ is the first projection. The boundary part $\partial_{r}\Omega$
is compact and equal to the manifold $\mathcal{E}_{r}$ of maximal tripotents.
\end{proposition}

The boundary part $\partial_{r}\Omega=\mathcal{E}_{r}$ is actually the \emph{%
Shilov boundary} of $\Omega$, that is, the smallest set of points where the
functions that are holomorphic on the domain and continuous up to the
boundary take their maximum modulus values.

The (affine) submanifold $\overline{e+\Omega(e)}$ is called \emph{affine
component} of $\partial\Omega$ (through the minimal tripotent $e$). It can
be shown that $\overline{e+\Omega(e)}$ is the maximal affine subset of $%
\partial\Omega$ containing $e$, which justifies its name. The decomposition (%
\ref{strat}) will be referred to as the \emph{stratification} of the
boundary. The submanifolds $\partial_{k}\Omega$ are called the \emph{%
boundary parts} (or \emph{strata}) of the boundary $\partial\Omega$.
Clearly, one has $\overline{\partial_{1}\Omega}=\partial\Omega$; the
submanifold $\partial _{1}\Omega$ has real codimension $1$ and is referred
to as the \emph{smooth part} of the boundary.

\subsubsection{The boundary of the exceptional domain of type VI}

Using Propositions \ref{prop3/7}, \ref{prop3/8} and the results about
tripotents from the previous section, we work out the details for the
exceptional Jordan triple system $\mathcal{H}_{3}(\mathbb{O})$.

\begin{proposition}
\label{prop3/9}Let $V=H_{3}(\mathbb{O)}$ be the exceptional Jordan triple
system of type $VI$. Then

1.~The manifold $\mathcal{E}_{1}$ of minimal tripotents is
\begin{align*}
\mathcal{E}_{1} & =\left\{ e\in V\mid\{e,e,e\}=2e,\ (e\mid e)=1\right\} \\
& =\left\{ e\in V\mid e^{\#}=0,\ (e\mid e)=1\right\} .
\end{align*}
For $e\in\mathcal{E}_{1}$, we have $V_{2}(e)=\mathbb{C}e$, $\dim_{\mathbb{C}%
}V_{1}(e)=16$. The manifold $\mathcal{E}_{1}$ has real dimension $33$ and is
a Cauchy-Riemann manifold of CR type $\left( 16,1\right) $.

2.~The manifold $\mathcal{E}_{2}$ of rank $2$ tripotents is defined by
\begin{align*}
\mathcal{E}_{2} & =\left\{ e\in V\mid\{e,e,e\}=2e,\ (e\mid e)=2\right\} \\
& =\left\{ e\in V\mid\det e=0,\ (e^{\#}\mid e^{\#})=1,\ (e\mid e)=2\right\} .
\end{align*}
For $e\in\mathcal{E}_{2}$, we have $\dim V_{2}(e)=10$ and $\dim V_{1}(e)=16$%
. The manifold $\mathcal{E}_{2}$ has real dimension $42$ and is a
Cauchy-Riemann manifold of CR type $\left( 16,10\right) $.

3.~The manifold $\mathcal{E}_{3}$ of maximal tripotents is defined by
\begin{align*}
\mathcal{E}_{3} & =\left\{ e\in V\mid\{e,e,e\}=2e,\ (e\mid e)=3\right\} \\
& =\left\{ e\in V\mid\left\vert \det e\right\vert ^{2}=1,\ (e^{\#}\mid
e^{\#})=3,\ (e\mid e)=3\right\} .
\end{align*}
The manifold $\mathcal{E}_{3}$ is totally real of real dimension $27$.
\end{proposition}

\begin{proof}
The characterization of the manifolds $\mathcal{E}_{k}$ has been obtained in
Proposition \ref{th2/10}. To study the spaces $V_{1}(e)$ and $V_{2}(e)$, we
use the fact that the group $K$ of automorphisms of the Jordan triple system
acts transitively on each $\mathcal{E}_{k}$ and that $u(V_{j}(e))=V_{j}(ue)$
for $u\in K$.

1.~Consider the minimal tripotent $e=e_{1}$. Then%
\begin{equation*}
V_{1}(e_{1})=\mathcal{F}_{2}\oplus\mathcal{F}_{3},\qquad V_{2}(e_{1})=%
\mathbb{C}e_{1},
\end{equation*}
which yields for all $e\in\mathcal{E}_{1}$, $V_{2}(e)=\mathbb{C}e$, $\dim_{%
\mathbb{C}}V_{1}(e)=16$.

2.~Consider the tripotent of rank $2$: $e=e_{1}+e_{2}$. As $e_{1}$ and $e_{2}
$ are orthogonal tripotents, we have
\begin{equation*}
D(e,e)=D(e_{1},e_{1})+D(e_{2},e_{2}).
\end{equation*}
From
\begin{align*}
V_{0}(e_{1}) & =\mathbb{C}e_{2}\oplus\mathbb{C}e_{3}\oplus\mathcal{F}_{1}, \\
V_{1}(e_{1}) & =\mathcal{F}_{2}\oplus\mathcal{F}_{3},\quad V_{2}(e_{1})=%
\mathbb{C}e_{1}
\end{align*}
and the analogous statement for $V_{j}(e_{2})$, we deduce
\begin{align*}
V_{0}(e) & =\mathbb{C}e_{3},\quad V_{1}(e)=\mathcal{F}_{1}\oplus \mathcal{F}%
_{2}, \\
V_{2}(e) & =\mathbb{C}e_{1}\oplus\mathbb{C}e_{2}\oplus\mathcal{F}_{3}.
\end{align*}
Hence we have $\dim V_{2}(e)=10$ and $\dim V_{1}(e)=16$ for all $e\in
\mathcal{E}_{2}$.

3.~Consider the maximal tripotent $e=e_{1}+e_{2}+e_{3}$. Then $V_{2}(e)=V$, $%
V_{0}(e)=V_{1}(e)=0.$ The tangent space direction to $\mathcal{E}_{3}$ at $e$
is $V_{2}^{+}(e)$ and is totally real of real dimension $27$.
\end{proof}

We now work out the specific information for the application of Proposition %
\ref{prop3/8} to the exceptional symmetric domain of dimension $27$.

\begin{proposition}
\label{prop3/3}1.\ The smooth boundary part $\partial_{1}\Omega$ is a
locally closed submanifold of real codimension $1$, defined by
\begin{align*}
f_{1}(x,x) & \equiv1-(x\mid x)+(x^{\#}\mid x^{\#})-\left\vert \det
x\right\vert ^{2}=0, \\
f_{2}(x,x) & \equiv3-2(x\mid x)+(x^{\#}\mid x^{\#})>0, \\
f_{3}(x,x) & \equiv3-(x\mid x)>0.
\end{align*}
It contains the manifold
\begin{equation*}
\mathcal{E}_{1}=\left\{ e\in H_{3}(\mathbb{O});(e\mid e)=1,\
e^{\#}=0\right\}
\end{equation*}
of minimal tripotents of $H_{3}(\mathbb{O})$. For $e\in\mathcal{E}_{1}$, the
normal direction at $e$ to $\partial_{1}\Omega$ is $e$ and the direction $%
\overrightarrow{T_{e}(\partial_{1}\Omega)}$ of the tangent space $%
T_{e}(\partial_{1}\Omega)$ is
\begin{equation*}
\overrightarrow{T_{e}(\partial_{1}\Omega)}=\operatorname{i}\mathbb{R}e\oplus
V_{1}(e)\oplus V_{0}(e).
\end{equation*}
The Peirce subspaces $V_{1}(e)$ and $V_{0}(e)$ have respective complex
dimensions $16$ and $10$. The Jordan triple subsystem $V_{0}(e)$ is
isomorphic to the classical Hermitian JTS of type $IV_{10}$ and the domain $%
\Omega (e)\subset V_{0}(e)$ is isomorphic to a Lie ball of dimension $10$.

2.~The boundary part $\partial_{2}\Omega$ is a locally closed,
Cauchy-Riemann submanifold of dimension $44$ and CR type $\left(
17,10\right) $; it contains the manifold
\begin{equation*}
\mathcal{E}_{2}=\left\{ e\mid(e\mid e)=2,\ (e^{\#}\mid e^{\#})=1,\ \det
e=0\right\}
\end{equation*}
of rank $2$ tripotents in $H_{3}(\mathbb{O})$. For $e\in\mathcal{E}_{2}$,
the normal direction $V_{2}^{+}(e)$ to $\partial_{2}\Omega$ at $e$ has real
dimension $10$; the direction $\overrightarrow{T_{e}(\partial_{2}\Omega)}$
of the tangent space $T_{e}(\partial_{2}\Omega)$ is
\begin{equation*}
\overrightarrow{T_{e}(\partial_{1}\Omega)}=\operatorname{i}V_{2}^{+}(e)\oplus
V_{1}(e)\oplus V_{0}(e),
\end{equation*}
where the Peirce subspaces $V_{1}(e)$ and $V_{0}(e)$ have respective complex
dimensions $16$ and $1$. The intersection of $\partial_{2}\Omega$ with the
affine tangent space $e+\overrightarrow{T_{e}(\partial_{2}\Omega)}$ is
\begin{equation*}
\partial_{2}\Omega\cap\left( e+\overrightarrow{T_{e}(\partial_{2}\Omega )}%
\right) =\partial_{2}\Omega\cap\left( e+V_{0}(e)\right) =e+\Omega(e),
\end{equation*}
where $\Omega(e)$ is the unit disc of the one dimensional Jordan triple
subsystem $V_{0}(e)$.

3.~The submanifold $\partial_{3}\Omega=\mathcal{E}_{3}$ is compact and
totally real (of real dimension $27$).
\end{proposition}

\begin{proof}
1.~Consider the minimal tripotent $e=e_{1}$. Then%
\begin{equation*}
V_{2}(e_{1})=\mathbb{C}e_{1},\quad V_{1}(e_{1})=\mathcal{F}_{2}\oplus
\mathcal{F}_{3},\quad V_{0}(e_{1})=\mathbb{C}e_{2}\oplus\mathbb{C}e_{3}\oplus%
\mathcal{F}_{1}.
\end{equation*}
In the Jordan triple subsystem $V_{0}(e_{1})$, a frame is given by $%
(e_{2},e_{3})$; the spaces of the Peirce decomposition of $T=V_{0}(e_{1})$
with respect to this frame are
\begin{align*}
& T_{2}(e_{2})  =\mathbb{C}e_{2},\quad T_{2}(e_{3})=\mathbb{C}e_{3},\quad
T_{1}(e_{2})\cap T_{1}(e_{3})=\mathcal{F}_{1}, \\
& T_{1}(e_{2})\cap T_{0}(e_{3})  =T_{0}(e_{2})\cap T_{1}(e_{3})=0.
\end{align*}
This shows that the Hermitian positive JTS is simple, with rank $r=2$ and
multiplicities $a=8$, $b=0$. The only possibility shown by the
classification of Hermitian positive JTS is the type $IV_{10}$. The
isomorphism of $V_{0}(e_{1})$ with the standard JTS of type $IV_{10}$ can
also be checked directly (see Exercise \ref{exerc3/1}).

2.~Consider the rank $2$ tripotent $e=e_{1}+e_{2}$. Then
\begin{equation*}
V_{0}(e) =\mathbb{C}e_{3},\quad V_{1}(e)=\mathcal{F}_{1}\oplus \mathcal{F}_{2},
\quad V_{2}(e) =\mathbb{C}e_{1}\oplus\mathbb{C}e_{2}\oplus\mathcal{F}_{3},
\end{equation*}
$\dim V_{2}(e)=10$ and $\dim V_{1}(e)=16$. The submanifold $%
\partial_{2}\Omega$ has normal direction $V_{2}^{+}(e)$, hence codimension $%
10$ and dimension $44$; the complex tangent direction to $\partial_{2}\Omega$
at $e$ is $V_{0}(e)\oplus V_{1}(e)=\mathbb{C}e_{3}\oplus\mathcal{F}_{1}\oplus%
\mathcal{F}_{2}$ and has dimension $17$, which means $\partial_{2}\Omega$
has CR type $\left( 17,10\right) $. The one-dimensional subsystem $V_{0}(e)=%
\mathbb{C}e_{3}$ admits as tripotent $e_{3}$, and we have $\Omega(e)=\Delta
e_{3}$, where $\Delta$ is the unit disc of $\mathbb{C}.$

3.~The equality $\partial_{3}\Omega=\mathcal{E}_{3}$ results from
Proposition \ref{prop3/8} and is also easily checked directly. The
properties of $\mathcal{E}_{3}$ have been obtained in Proposition \ref%
{prop3/9}.
\end{proof}

\begin{exercise}
\label{exerc3/1}Compute the Jordan triple product in $V_{0}(e_{1})$ and show
that the JTS $V_{0}(e_{1})$ is isomorphic to the Hermitian positive JTS of
type $IV_{10}$.
\end{exercise}

\begin{exercise}
\label{exerc3/2}1.~For a minimal tripotent $e=\mathcal{E}_{1}$, compute
explicitly the Peirce subspaces $V_{1}(e)$ and $V_{0}(e)$ and their Jordan
triple structure.

2.~Compute the map
\begin{equation*}
p_{1}:\partial_{1}\Omega\rightarrow\mathcal{E}_{1}
\end{equation*}
defined by $p_{1}(x)=e$ if $x\in\mathcal{E}_{1}$ and $x-e\in V_{0}(e)$,
using the operations in $\mathcal{H}_{3}(\mathbb{O})$.
\end{exercise}

\begin{exercise}
\label{exerc3/3}1.~For a rank $2$ tripotent $e=\mathcal{E}_{2}$, compute
explicitly the Peirce subspaces $V_{1}(e)$ and $V_{0}(e)$ and their Jordan
triple structure.

2.~Compute the map
\begin{equation*}
p_{2}:\partial_{2}\Omega\rightarrow\mathcal{E}_{2}
\end{equation*}
defined by $p_{2}(x)=e$ if $x\in\mathcal{E}_{2}$ and $x-e\in V_{0}(e)$,
using the operations in $\mathcal{H}_{3}(\mathbb{O})$.
\end{exercise}

\subsubsection{The boundary of the exceptional domain of type V}

Along the same lines, we describe the manifolds of tripotents of the
exceptional Jordan triple system of type $V$ and the boundary of the
associated symmetric domain.

\begin{proposition}
\label{prop3/10}Let $W=V_{1}(e_{1})=\mathcal{F}_{2}\oplus\mathcal{F}%
_{3}\subset H_{3}(\mathbb{O)}$ be the exceptional Jordan triple system of
type $V$. Then

1.~The manifold $\mathcal{E}_{1}^{\prime}$ of minimal tripotents is defined
by
\begin{align*}
\mathcal{E}_{1}^{\prime} & =\left\{ e\in W\mid\{e,e,e\}=2e,\ (e\mid
e)=1\right\} \\
& =\left\{ e\in W\mid e^{\#}=0,\ (e\mid e)=1\right\} .
\end{align*}
For $e\in\mathcal{E}_{1}^{\prime}$, we have $W_{2}(e)=\mathbb{C}e$, $\dim_{%
\mathbb{C}}W_{1}(e)=10$. The manifold $\mathcal{E}_{1}^{\prime}$ has real
dimension $21$ and is a Cauchy-Riemann manifold of CR type $\left(
10,1\right) $.

2.~The manifold $\mathcal{E}_{2}^{\prime}$ of maximal tripotents is defined
by
\begin{align*}
\mathcal{E}_{2}^{\prime} & =\left\{ e\in W\mid\{e,e,e\}=2e,\ (e\mid
e)=2\right\} \\
& =\left\{ e\in W\mid(e^{\#}\mid e^{\#})=1,\ (e\mid e)=2\right\} .
\end{align*}
For $e\in\mathcal{E}_{2}^{\prime}$, we have $\dim W_{2}(e)=8$ and $\dim
W_{1}(e)=8$. The manifold $\mathcal{E}_{2}^{\prime}$ has real dimension $24$
and is a Cauchy-Riemann manifold of CR type $\left( 8,8\right) $.
\end{proposition}

\begin{proof}
The description of $\mathcal{E}_{1}^{\prime}$ and $\mathcal{E}_{2}^{\prime}$
has been given in Proposition \ref{th2/18}.

1.~From Lemma \ref{th2/16}, we can take as minimal tripotent $u=F_{2}(\beta )
$, where $\beta\in\mathbb{O}$ such that $(\beta\mid\beta)=1$ and $n(\beta )=0
$. The spaces of the Peirce decomposition of $W$ with respect to $u$ are
\begin{align*}
W_{0}(u) & =\mathbb{C}\overline{u}\oplus F_{3}(\ker L\left( \overline {\beta}%
\right) ), \\
W_{1}(u) & =F_{2}(\left\langle \beta,\overline{\beta}\right\rangle ^{\perp
})\oplus F_{3}(\ker L\left( \beta\right) ), \\
W_{2}(u) & =\mathbb{C}u.
\end{align*}
From Lemma \ref{th2/15}, we know that
\begin{equation*}
\dim_{\mathbb{C}}\ker L(\beta)=\dim_{\mathbb{C}}\ker L(\overline{\beta})=4,
\end{equation*}
which yields $\dim W_{0}(u)=5$ and $\dim W_{1}(u)=10$.

2.~We can choose as maximal tripotent $w=u+v=F_{2}(\beta+\overline{\beta }%
)=F_{2}(c)$, with $c\in\mathbb{O}_{c}$ and $n(c)=1$. The Peirce spaces for $w
$ are $W_{2}(w)=\mathcal{F}_{2}$, $W_{1}(w)=\mathcal{F}_{3}$; both have
complex dimension $8$.
\end{proof}

\begin{proposition}
\label{prop3/2}The boundary $\partial\Omega^{\prime}$ of the exceptional
symmetric domain $\Omega^{\prime}=\Omega_{W}$ of type $V$ is the disjoint
union%
\begin{equation}
\partial\Omega^{\prime}=\partial_{1}\Omega^{\prime}\amalg\partial_{2}%
\Omega^{\prime}   \label{strat-16}
\end{equation}
of the locally closed manifold $\partial_{1}\Omega^{\prime}$ and of the
compact manifold $\partial_{2}\Omega^{\prime}$.

1.~The smooth boundary part $\partial_{1}\Omega^{\prime}$ is a locally
closed submanifold of real codimension $1$; it contains the manifold%
\begin{equation*}
\mathcal{E}_{1}^{\prime}=\left\{ x\in W;\ (x\mid x)=1,\ x^{\#}=0\right\} .
\end{equation*}
of minimal tripotents of $W=\mathcal{F}_{2}\oplus\mathcal{F}_{3}$. For $e\in%
\mathcal{E}_{1}^{\prime}$, the normal direction at $e$ to $%
\partial_{1}\Omega^{\prime}$ is $e$; the direction $\overrightarrow {%
T_{e}(\partial_{1}\Omega^{\prime})}$ of the tangent space $T_{e}(\partial
_{1}\Omega^{\prime})$ is
\begin{equation*}
\overrightarrow{T_{e}(\partial_{1}\Omega^{\prime})}=\operatorname{i}\mathbb{R}%
e\oplus W_{1}(e)\oplus W_{0}(e),
\end{equation*}
where the Peirce subspaces $W_{1}(e)$ and $W_{0}(e)$ have respective complex
dimensions $10$ and $5$. The affine component $e+\Omega^{\prime}(e)$ is the
unit Hermitian ball with center $e$ in $e+W_{0}(e)$.

2.~The submanifold $\partial_{2}\Omega^{\prime}=\mathcal{E}_{2}^{\prime}$ is
a compact, Cauchy-Riemann submanifold of CR type $(16,8)$ and real dimension
$24$.
\end{proposition}

\begin{proof}
1.~ Consider the minimal tripotent $u=F_{2}(\beta)$, where $\beta\in \mathbb{%
O}$ such that $(\beta\mid\beta)=1$ and $n(\beta)=0$. From Lemma \ref{th2/16}%
, we know that the spaces of the Peirce decomposition of $W$ with respect to
$u$ are
\begin{align*}
W_{0}(u) & =\mathbb{C}\overline{u}\oplus F_{3}(\ker L\left( \overline {\beta}%
\right) ), \\
W_{1}(u) & =F_{2}(\left\langle \beta,\overline{\beta}\right\rangle ^{\perp
})\oplus F_{3}(\ker L\left( \beta\right) ), \\
W_{2}(u) & =\mathbb{C}u.
\end{align*}
This proves the statement about the tangent space and the dimensions of
Peirce subspaces.

Let $v$ be a tripotent in $W_{0}(u)$. Then $v$ is orthogonal to $u$ and $%
(u,v)$ is a frame of $W$. This implies that $v$ is maximal in $W_{0}(u)$ and
that the JTS $W_{0}(u)$ is of rank $1$, so that $\Omega(u)$ is a Hermitian
ball in $W_{0}(u)$.

2.~An element $x\in W$ belongs to $\partial_{2}\Omega^{\prime}$ if and only
if
\begin{align*}
g_{1}(x,x) & \equiv1-(x\mid x)+(x^{\#}\mid x^{\#})=0, \\
g_{2}(x,x) & \equiv2-(x\mid x)=0.
\end{align*}
These conditions are clearly equivalent to%
\begin{equation*}
(x\mid x)=2,\quad(x^{\#}\mid x^{\#})=1,
\end{equation*}
which is precisely the characterization of elements of $\mathcal{E}%
_{2}^{\prime}$. So $\partial_{2}\Omega^{\prime}=\mathcal{E}_{2}^{\prime}$,
which also results from the general theory. The structure of $\mathcal{E}%
_{1}^{\prime}$ has been given in Proposition \ref{prop3/10}.
\end{proof}

\begin{exercise}
\label{exerc3/5}1.~For a minimal tripotent $e=F_{2}(\beta)+F_{3}(\gamma )\in%
\mathcal{E}_{1}^{\prime}$, compute explicitly the Peirce subspaces $W_{1}(e)$
and $W_{0}(e)$ and their Jordan triple structure.

2.~Compute explicitly the map
\begin{equation*}
p_{1}^{\prime}:\partial_{1}\Omega^{\prime}\rightarrow\mathcal{E}_{1}^{\prime}
\end{equation*}
defined by $p_{1}^{\prime}(x)=e$ if $x\in\mathcal{E}_{1}^{\prime}$ and $%
x-e\in W_{0}(e)$.

3.~Identify the type of the Hermitian positive JTS $W_{1}(e)$, where $e$ is
a minimal tripotent of $W$.
\end{exercise}

\begin{exercise}
\label{exerc3/4}For a maximal tripotent $e=F_{2}(\beta)+F_{3}(\gamma )\in%
\mathcal{E}_{2}^{\prime}$, compute explicitly the Peirce subspaces $W_{2}(e)$
and $W_{1}(e)$ and their Jordan triple structure.
\end{exercise}

\subsection{Compactification of exceptional symmetric domains}

\index{compactification}In this section, we work out the canonical
projective realization of the compact dual of the two exceptional domains
(see \cite{Loos1977}, \cite{Roos1999}).

\subsubsection{The Freudenthal manifold}

Consider the exceptional Jordan triple $V=H_{3}(\mathbb{O)}$ of type $VI$.
The \emph{generic norm} of $V$ is
\begin{equation*}
N_{V}(x,y)=m_{V}(1,x,y)=1-(x\mid y)+(x^{\#}\mid y^{\#})-\det x\det%
\overline {y}.
\end{equation*}
To this generic norm, we associate the map
\begin{align*}
j :V & \rightarrow\mathbb{P}\left( \mathbb{C\oplus}V\mathbb{\oplus }V\mathbb{%
\oplus C}\right) \\*
x & \mapsto\left[ 1,x,x^{\#},\det x\right] ,
\end{align*}
where $\left[ \ldots\right] $ denotes the class in the projective space.

\begin{definition}
The \index{FreudenthalH@\textsc{Freudenthal, Hans} (1905--1990)}%
\index{Freudenthal manifold}\emph{Freudenthal manifold} is the submanifold
of $\,\mathbb{P}(\mathbb{C\oplus}V\mathbb{\oplus}V\mathbb{\oplus C})$ defined by%
\begin{equation*}
\mathcal{M}=\{ \left[ \lambda,x,y,\mu\right] \mid\lambda,\mu \in\mathbb{C},\
x,y\in V,\ y^{\#}=\mu x,\ x^{\#}=\lambda y,\ (x:y)=3\lambda \mu\} .
\end{equation*}
\end{definition}

Note that this definition makes sense, for the defining equations $%
y^{\#}=\mu x$, $x^{\#}=\lambda y$ and $(x:y)=3\lambda\mu$ are homogeneous of
degree $2$. As $( x^{\#}) ^{\#}=x\det x$ and $\left(
x:x^{\#}\right) =3\det x$ for $x\in V$, we have $j(V)\subset\mathcal{M}$.
The map $j$ is clearly an immersion.

Let $\left[ \lambda,x,y,\mu\right] \in\mathcal{M}$ and assume $\lambda\neq 0$%
. Let $x^{\prime}= \displaystyle\frac{x}{\lambda}$; then%
\begin{align*}
\left( x^{\prime}\right) ^{\#} & =\frac{x^{\#}}{\lambda^{2}}=\frac {y}{%
\lambda}, \\
\det(x^{\prime}) & =\frac{\det x}{\lambda^{3}}=\frac{1}{3}\frac{(x^{\#}:x)}{%
\lambda^{3}}=\frac{1}{3}\frac{(y:x)}{\lambda^{2}}=\frac{\mu}{\lambda},
\end{align*}
which shows that $\left[ \lambda,x,y,\mu\right] =j(x^{\prime})$.

\begin{proposition}
\label{compact-27}The map $j$ is an immersion of $\,V$ onto an open dense
subset of the Freudenthal manifold $\mathcal{M}$.
\end{proposition}

\subsubsection{Compactification of the $16$-dimensional exceptional domain}

Consider the exceptional symmetric domain of dimension $16$, realized as%
\begin{equation*}
W=\mathcal{F}_{2}\oplus\mathcal{F}_{3}\subset V=H_{3}(\mathbb{O}).
\end{equation*}
The \emph{generic norm} of $W$ is
\begin{equation*}
N_{W}(x,y)=m_{W}(1,x,y)=1-(x\mid y)+(x^{\#}\mid y^{\#}).
\end{equation*}
One checks easily from the definition of $x^{\#}$ that $x\in W$ implies $%
x^{\#}\in V_{0}(e_{1})$. Note that the Peirce decomposition of $V$ with
respect to $e_{1}$ has the eigenspaces
\begin{align*}
V_{2}(e_{1}) & =\mathbb{C}e_{1},\ V_{1}(e_{1})=W=\mathcal{F}_{2}\oplus%
\mathcal{F}_{3}, \\
V_{0}(e_{1}) & =\mathbb{C}e_{2}\oplus\mathbb{C}e_{3}\oplus\mathcal{F}_{1}.
\end{align*}

\begin{lemma}
Let $z=e_{1}+x+y$ with $x\in W$ and $y\in V_{0}(e_{1})$. Then $z^{\#}=0$ if
and only if $y=-e_{1}\times x^{\#}$.
\end{lemma}

\begin{proof}
Let $z=e_{1}+x+y$, $x=F_{2}(b)+F_{3}(c)\in W$, $y=\mu e_{2}+\nu
e_{3}+F_{1}(a)\in V_{0}(e_{1})$. We have
\begin{align*}
e_{1}^{\#} & =0,\ x^{\#}=-n(b)e_{2}-n(c)e_{3}+\widetilde{F_{1}}(bc), \\
y^{\#} & =\left( \mu\nu-n(a)\right) e_{1},\ e_{1}\times x=0, \\
e_{1}\times y & =\mu e_{3}+\nu e_{2}-F_{1}(a), \\
x\times y & =-\mu F_{2}(b)-\nu F_{3}(c)+\widetilde{F_{3}}(ab)+\widetilde {%
F_{2}}(ca)
\end{align*}
and%
\begin{align*}
z^{\#} & =\left( e_{1}+x+y\right) ^{\#}=e_{1}^{\#}+x^{\#}+y^{\#}+e_{1}\times
x+e_{1}\times y+x\times y \\*
& =-n(b)e_{2}-n(c)e_{3}+\widetilde{F_{1}}(bc)+\left( \mu\nu-n(a)\right) e_{1}
\\*
& \quad +\mu e_{3}+\nu e_{2}-F_{1}(a)-\mu F_{2}(b)-\nu F_{3}(c)+\widetilde{F_{3}}%
(ab)+\widetilde{F_{2}}(ca) \\
& =\left( \mu\nu-n(a)\right) e_{1}+\left( \nu-n(b)\right) e_{2}+\left(
\mu-n(c)\right) e_{3} \\*
& \quad +\widetilde{F_{1}}(bc)-F_{1}(a)+\widetilde{F_{2}}(ca)-\mu F_{2}(b)+%
\widetilde{F_{3}}(ab)-\nu F_{3}(c).
\end{align*}

Then $z^{\#}=0$ implies%
\begin{equation}
\mu=n(c),\ \nu=n(b),\ a=\widetilde{bc},   \label{eq3/4}
\end{equation}
that is,
\begin{equation*}
y=n(c)e_{2}+n(b)e_{3}+\widetilde{F_{1}}(bc),
\end{equation*}
which is equivalent to
\begin{equation*}
y=-e_{1}\times x^{\#}.
\end{equation*}

Conversely, if $y=-e_{1}\times x^{\#}$, the relations (\ref{eq3/4}) are
satisfied and imply
\begin{align*}
&n(a) =n(bc)=n(b)n(c)=\mu\nu, \\
& \widetilde{ca}  =(bc)\widetilde{c}=\mu b, \quad
\widetilde{ab}  =\widetilde{b}(bc)=\nu c,
\end{align*}
which shows that $z^{\#}=0$.
\end{proof}

With the help of the previous lemma, we are now able to describe a
compactification of $W$, isomorphic to the canonical compactification
associated to the generic norm $N_{W}$.

\begin{proposition}
Let $V=H_{3}(\mathbb{O)}$ and $W=\mathcal{F}_{2}\oplus\mathcal{F}_{3}$.
Define $j:W\rightarrow\mathbb{P}(V)$ by
\begin{equation*}
j(x)=\left[ e_{1}+x-e_{1}\times x^{\#}\right] \quad(x\in W).
\end{equation*}
Then $j$ is a biholomorphism of $W$ onto an open dense subset of the manifold%
\begin{equation*}
\mathcal{P}=\left\{ [z]\in\mathbb{P}(V)\mid z^{\#}=0\right\} .
\end{equation*}
\end{proposition}

Indeed, $j$ is an immersion and maps $W$ biholomorphically onto
\begin{equation*}
j(W)=\left\{ [z]\in\mathcal{P}\mid(z:e_{1})\neq0\right\}.
\end{equation*}

The manifold $\mathcal{P}$ is the image in $\mathbb{P}(V)$ of the cone $%
\left\{ z^{\#}=0\right\} $ of rank one elements in $V$.


\begin{theindex}

  \item $\mathbb{H}$, 11
  \item $\mathcal{H}_{3}(\mathbb{O}_{\mathbb{C}})$, 2, 12
  \item $\mathcal{M}_{2,1}(\mathbb{O}_{\mathbb{C}})$, 2
  \item $\mathbb{O}$, $\mathbb{O}_{c}$, 11
  \item $\mathbb{O}_{\mathbb{C}}$, 11
  \item $\mathbb{O}_{s}$, 11
  \item $a^{\#}$, adjoint of $a\in H_{3}(\mathbb{O})$, 12

  \indexspace

  \item Albert algebra, 2
  \item \textsc{Albert, Abraham Adrian} (1905--1972), 2
  \item alternative algebra, 5
  \item alternativity, 5
  \item associator, 5

  \indexspace

  \item Cayley algebra
    \subitem compact ---, 11
    \subitem complex ---, 11
    \subitem split ---, 11
  \item Cayley conjugation, 4
  \item Cayley--Dickson extension, 7
  \item \textsc{Cayley, Arthur} (1821--1895), 2, 7
  \item compactification, 32
  \item composition algebra, 3
    \subitem opposite ---, 3

  \indexspace

  \item determinant in $H_{3}(\mathbb{O})$, 13
  \item \textsc{Dickson, Leonard Eugene} (1874--1954), 7

  \indexspace

  \item exceptional Jordan algebra, 2
  \item exceptional JTS
    \subitem of dimension 16, 24
    \subitem of dimension 27, 17
  \item exceptional symmetric domain
    \subitem of dimension 16, 26
    \subitem of dimension 27, 26

  \indexspace

  \item flat subspace in Hermitian JTS, 17
  \item flexible algebra, 6
  \item Freudenthal manifold, 32
  \item Freudenthal product, 12
  \item Freudenthal's theorem, 21
  \item \textsc{Freudenthal, Hans} (1905--1990), 12, 21, 32

  \indexspace

  \item generic minimal polynomial, 17
  \item \textsc{Graves, John Thomas} (1806--1870), 2

  \indexspace

  \item \textsc{Hua, Luokeng} (1910--1985), 2
  \item Hurwitz algebra, 3
  \item \textsc{Hurwitz, Adolf} (1859--1919), 3

  \indexspace

  \item Jordan algebra, 2
  \item Jordan triple product, 2
  \item Jordan triple system, 2, 15
  \item \textsc{Jordan, Pascual} (1902--1980), 2

  \indexspace

  \item Moufang identities, 6
  \item \textsc{Moufang, Ruth} (1905--1977), 6

  \indexspace

  \item norm
    \subitem  on a composition algebra, 3

  \indexspace

  \item octonions, 11
    \subitem complex ---, 11
  \item odd powers in Hermitian JTS, 17
  \item orthogonal tripotents, 18

  \indexspace

  \item Peirce decomposition, 20, 26
  \item \textsc{Peirce, Charles Sanders} (1839--1914), 20

  \indexspace

  \item rank of a Jordan triple, 17

  \indexspace

  \item tripotent in Hermitian JTS, 17

\end{theindex}


\end{document}